\def\version{24.2.2022}
\def\users{us}  %
\def\users{final-layout}   % when activated, ``our'' debugging is suppressed
\documentclass[12pt]{article}
\usepackage{upgreek}
\usepackage{xcolor}
\usepackage{bm,amsmath,amsthm,hyperref,amsfonts,amssymb,color}
\usepackage{mathrsfs} % THIS PACKAGE ALLOWS FOR ``scr'' FONTS
\usepackage{cite}

\usepackage{ifthen}
\ifthenelse{\equal{\users}{final-layout}}{}{
\usepackage{fancyhdr}
\pagestyle{fancy}
\headheight=28pt\headwidth=17cm
\definecolor{gray}{gray}{0.5}
\rhead{\color{gray}Quasistatic hypoplasticity\\
T.Roub\'\i\v cek }
\chead{}
\lhead{Version \version, file: \jobname.tex,
\\
compiled:
\number\day.\number\month.\number\year\ at
\the\hour:\ifnum\minute<10 0\fi\the\minute\ h }
}

\definecolor{labelkey}{rgb}{1.,.2,0.}
\newcount\hour \newcount\minute
\hour=\time
\divide \hour by 60
\minute=\time
\loop \ifnum \minute > 59 \advance \minute by -60 \repeat

\usepackage[normalem]{ulem}

\usepackage[normalem]{ulem}
\usepackage{ifthen}
\usepackage{color}

\ifthenelse{\equal{\users}{final-layout}}{

	\newcommand{\COMMENT}[1]{}
	\newcommand{\COMMENTGT}[1]{}
	\newcommand{\TODO}[1]{}
	\newcommand{\INTERNAL}[1]{}
	\newcommand{\QUESTION}[1]{}
	\newcommand{\DELETE}[1]{}

	\newcommand{\REM}[1]{\marginpar{\bfseries\tiny{\color{blue}}}}
    \newcommand{\MARGINOTE}[1]{}
}
{
	
	\newcommand{\COMMENT}[1]{{\color{red}\uuline{#1}\color{black}}}
	\newcommand{\COMMENTGT}[1]{{\hfill\large\color{red}***{#1}***\color{black}\hfill}\\}
	\newcommand{\TODO}[1]{{\color{red}\uuline{#1}\color{black}}}
	\newcommand{\INTERNAL}[1]{\footnote{#1}}
	\newcommand{\QUESTION}[1]{{\color{brown}\uuline{#1}\color{black}}}
	\newcommand{\DELETE}[1]{{\color{red}\sout{#1}\color{black}}}

	\newcommand{\REM}[1]{\marginpar{\bfseries\tiny{\color{blue}#1}}}
\newcommand{\MARGINOTE}[1]{\marginpar{\color{red}\tiny\texttt{#1}}}
}

\newcommand\DT[1]{\mathchoice
                 {{\buildrel{\hspace*{.1em}\text{\LARGE.}}\over{#1}}}
                 {{\buildrel{\hspace*{.1em}\text{\Large.}}\over{#1}}}
                 {{\buildrel{\hspace*{.1em}\text{\large.}}\over{#1}}}
                 {{\buildrel{\hspace*{.1em}\text{\large.}}\over{#1}}}}
\newcommand\pdt[1]{\frac{\partial{#1}}{\partial t}} %Partial Derivative w.r.t. t
\newcommand{\lineunder}[2]{\LU{\begin{array}[t]{c}\underbrace{#1}\vspace*{.5em}\end{array}}{\mbox{\footnotesize\rm #2}}}
\newcommand{\LU}[2]{\begin{array}[t]{c}#1\vspace*{-1em}\\_{#2}\end{array}}
\newcommand{\linesunder}[3]{\LSU{\begin{array}[t]{c}\underbrace{#1}\vspace*{.5em}\end{array}}{\mbox{\footnotesize\rm #2}}{\mbox{\footnotesize\rm #3}}}
\newcommand{\LSU}[3]{\begin{array}[t]{c}#1\vspace*{-1em}\\_{#2}\vspace*{-.5em}\\_{#3}\end{array}}
\newcommand{\morelinesunder}[4]{\LSUU{\begin{array}[t]{c}\underbrace{#1}\vspace*{.5em}\end{array}}{\mbox{\footnotesize\rm #2}}{\mbox{\footnotesize\rm #3}}{\mbox{\footnotesize\rm #4}}}
\newcommand{\LSUU}[4]{\begin{array}[t]{c}#1\vspace*{-1em}\\_{#2}\vspace*{-.5em}\\_{#3}\vspace*{-.5em}\\_{#4}\end{array}}
\newcommand{\Item}[2]{\parbox[t]{.055\textwidth}{#1}\hfill%
      \parbox[t]{.945\textwidth}{#2}\vspace*{.8mm}} 
\newcommand{\divS}{\mathrm{div}_{\scriptscriptstyle\textrm{\hspace*{-.1em}S}}^{}}
\newcommand{\nablaS}{\nabla_{\scriptscriptstyle\textrm{\hspace*{-.3em}S}}^{}}
\newcommand{\NablaS}{\Nabla_{\scriptscriptstyle\textrm{\hspace*{-.3em}S}}^{}}
\def\Vdots{\!\mbox{\setlength{\unitlength}{1em}
\begin{picture}(0,0)
\put(-.07,0){.}
\put(-.07,.3){.}
\put(-.07,.6){.}
\end{picture}\hspace*{.2em}}}
\usepackage{mathrsfs}   % loading \mathscr
\usepackage{eucal}      % changing \mathcal into EulerCaligraphic  

  \def\bbI{{\mathbb I}}

\def\FG{\boldsymbol}
   
 \def\ee{{\FG e}} \def\ff{{\FG f}}

 \def\nn{{\FG n}}

\def\vv{{\FG v}}  \def\xx{{\FG x}} 
\def\yy{{\FG y}}  

\def\FF{{\FG F}} 
  \def\LL{{\FG L}}

\def\SS{{\FG S}} \def\TT{{\FG T}} %\def\UU{{\FG U}} 
  \def\XX{{\FG X}}

\newcommand{\R}{\mathbb R}
\newcommand{\N}{\mathbb N}
\newcommand{\Nabla}{{\bm\nabla}}
\newcommand{\Fe}{\FF_{\hspace*{-.2em}\mathrm e^{^{^{}}}}}
\newcommand{\Fee}{\FF_{\hspace*{-.2em}\mathrm e}}
\newcommand{\Fekl}{\FF_{{\hspace*{-.2em}\mathrm e},kl}}
\newcommand{\FEk}{\FF_{{\hspace*{-.2em}\mathrm e},k}}

\newcommand{\Fezero}{\FF_{\hspace*{-.2em}\mathrm e,0}^{}}
\newcommand{\Fp}{\FF_{\hspace*{-.2em}\mathrm p^{^{^{}}}}}
\newcommand{\widetildeFekl}{\hspace*{.1em}\widetilde{\hspace*{-.1em}\FF}_{\hspace*{-.1em}{\mathrm e},kl^{^{^{}}}}}
\newcommand{\widetildeLp}{\widetilde{\LL\,}_{\hspace*{-.2em}\mathrm p^{^{^{}}}}}
\newcommand{\Fpzero}{\FF_{\hspace*{-.2em}\mathrm p^{^{^{}}},0}}
\newcommand{\Fpp}{\FF_{\hspace*{-.2em}\mathrm p}}

\newcommand{\pl}{\partial}
\newcommand{\eq}[1]{(\ref{#1})}
\newcommand{\Lp}{\LL_{\mathrm p^{^{^{}}}}}

\newcommand{\Lpk}{\LL_{\hspace*{-.0em}{\mathrm p},k}}
\newcommand{\Lpkl}{\LL_{\hspace*{-.0em}{\mathrm p},kl}}
\renewcommand{\d}{\mathrm d}  
\newcommand{\barOmega}{\hspace*{.2em}{\overline{\hspace*{-.2em}\varOmega}}}

\newcommand{\DTFe}{\DT\FF_{\!\mathrm e^{^{^{}}}}}
\newcommand{\DTFp}{\DT\FF_{\!\mathrm p^{^{^{}}}}}

\newtheorem{theorem}{Theorem}[section]
\newtheorem{lemma}[theorem]{Lemma}
\newtheorem{definition}[theorem]{Definition}

\newtheorem{proposition}[theorem]{Proposition}

\newtheorem{remark}[theorem]{Remark}

\numberwithin{equation}{section}
\def\NU{\nu}
\def\MU{\mu}

\def\GRAVITY{{\bm g}}

\begin{document}

\allowdisplaybreaks

\bigskip\bigskip\bigskip

\noindent{\Large\bf Quasistatic hypoplasticity
%\\[.2em]
at large strains Eulerian}

\bigskip\bigskip\bigskip\bigskip

\noindent{\bf Tom\'a\v s Roub\'\i\v cek}\footnote{Mathematical Institute, Charles University,
Sokolovsk\'a 83, CZ-186~75~Praha~8,  Czech Republic,\\\hspace*{1.6em} email: ${\texttt{tomas.roubicek@mff.cuni.cz}}$}
\footnote{Institute of Thermomechanics, Czech Academy of Sciences,
Dolej\v skova 5, CZ-18200~Praha~8, Czech Rep.}

%\end{center}

\bigskip\bigskip\bigskip

{\small
\baselineskip=12pt
\noindent{\bf Abstract}.
The isothermal quasistatic (i.e.\ acceleration neglected)
hardening-free plasticity at large strains is considered, based on the standard
multiplicative decomposition of the total strain and the isochoric plastic
distortion. The Eulerian velocity-strain formulation is used. The mass density
evolves too, but acts only via the force term with a given external acceleration.
This rather standard model is then re-formulated in terms of rates (so-called
hypoplasticity) and the plastic distortion is completely eliminated, although
it can be a-posteriori re-constructed. Involving gradient theories for dissipation,
existence and regularity of weak solutions is proved rather constructively
by a suitable regularization combined with a Galerkin approximation.
The local non-interpenetration through a blowup of stored energy when
elastic-strain determinant approaches zero is enforced and exploited. 
The plasticity is considered rate dependent and, as a special case, also creep
in Jeffreys' viscoelastic rheology in the shear is covered while the volumetric
response obeys the Kelvin-Voigt rheology.

\medskip


\medskip

\noindent {\it AMS Subject Classification:} 
35Q49, % Transport equations
35Q74, % PDEs in connection with mechanics of deformable solids
%35Q79, % PDEs in connection with classical thermodynamics and heat transfer
%35Q86, % PDEs in connection with geophysics
65M60, % PDE, IVP, ... Galerkin methods, 
%74A15, % Thermodynamics
74A30, % Nonsimple materials
% 74A45, % Theories of fracture and damage
74C15, % Large-strain, rate-independent theories(including nonlinear plasticity)
74Dxx. %Materials of strain-rate type.
%74F10, % Fluid-solid interactions (...aero- hydro-elasticity, porosity, etc.)
%74J30, % Nonlinear waves
%74L05, % Geophysical solid mechanics
%74R20, % Anelastic fracture and damage
%76S05, % Flows in porous media; filtration; seepage
%80A20. % Heat and mass transfer, heat flow
%86A17. % Global dynamics, earthquake problems

}

\bigskip\bigskip

\baselineskip=16pt

\section{Introduction}

Many materials undergo large inelastic process, specifically
plastification or creep. Typically, beside metals (some of whose
exhibit so-called ``superplasticity''), it concerns polymers and
particularly geomaterials as rocks, soils, and ice which can
exhibit very large inelastic strains on long time scales. Mechanically,
{\it large-strain} (sometimes called finite-strain) {\it plasticity}
or {\it creep} models have been developed during decades, see the monographs
\cite{Bert08EPLD,BesGie94MMID,GuFrAn10MTC,Hash20NCMF,HasYam13IFST,KhaHua95CTP,JirBaz02IAS,Luba02ET,Maug92TPF}
and references therein. Following a general idea to express
strain-stress responses rather in terms of rates, like hypo-elasticity
being an alternative description to hyperelasticity,
%rate-type materials (Rajagopal \& co.)
a rate formulation (sometimes called {\it hypoplasticity}) has been developed
as an alternative to the classical theory of
elasto-plasticity, cf.\ \cite{Dafa86PSCS,Koly91OH,LCCV04GFH},
 although the label ``hypoplasticity'' has rather free meaning
and is used in various ways, even not entirely identical as here, cf.\ also
%or \cite[Sect.24.2.4]{JirBaz02IAS}
\cite{TaViCh00RTDA}. 
One attribute is that this formulation
works without any explicit decomposition of the strain rate tensor to
a reversible and an irreversible parts, although it is implicitly based on it.

Let us summarize the main ingredients, which are actually quite
standard and generally accepted, and which will be employed:\\
\Item{---}{deformation in actual configuration (i.e.\ {\it Eulerian
approach}) and corresponding evolution of the deformation gradient, cf.\
\eq{ultimate} below,}
\Item{---}{corresponding transport of mass (i.e.\ {\it mass conservation}),
cf.\ \eq{cont-eq+},}
\Item{---}{Lie-Liu-Kr\"oner
%Green-Naghdi's
{\it multiplicative decomposition} of the
deformation gradient to the elastic and the inelastic (plastic) strains, cf.\ 
\eq{Green-Naghdi}, with the plastic distortion being {\it isochoric}, i.e.\
having determinant equal 1,}
\Item{---}{a stored energy dependent on elastic strain, expected generally
{\it nonconvex, frame indifferent, and singular} when respecting local
non-interpenetration by a blow-up within infinite compression, i.e.\ if
determinant of the elastic strain goes to zero, cf.\ \eq{Euler-ass-phi},}
\Item{---}{a dissipation potential acting on  the  symmetric velocity
gradient and on the plastic distor{{t}}ion rate,}
\Item{---}{the evolution based on the mentioned conservation of mass
and evolution of the deformation gradient, in addition on the momentum
equilibrium and on the flow rule of plastic distortion through the plastic
distor{{t}}ion rate, and}
\Item{---}{gradient theories, here applied on the dissipative potential.}
Ultimately, we focus to an energetics of the models, which will make a solid
base for %the
 a  rigorous analysis.

Any reference configuration (i.e.\ the Lagrangian approach as in
\cite{DaRoSt??NHFV,KruRou19MMCM,MaiMie09GERI,MieRou16RIEF,MiRoSa18GERV,RouSte19FTCS})
is thus eliminated from the formulation of the problem. This is very natural
especially for materials where such a reference configuration
cannot be identified naturally, as\ e.g.\ in geological materials
(rocks, soils, ice, etc.) which are permanently evolving on long time scales and which
do not possess any ``original'' stress-free configuration, cf.\
e.g.\ \cite{Niem03EHMS}. Rather, they have a continuously evolving
natural configuration, sometimes presented under the name of ``multiple natural
configurations'' \cite{RajSri04TMMN,RajSri05REEM}.
This is one of conventional approaches to inelasticity, dated back
to C.\,Eckart \cite{Ecka48TIPT}, including both creep and
plasticity. Rather for explanatory lucidity, we will present it in
detail in Section~\ref{sec-classic} first a classical way
including inelastic (plastic) distor{{t}}ion and the multiplicative
decomposition of the deformation gradient. The plastic isochoricity
is build in  through the dissipation potential.

Then, in Section~\ref{sec-hypo} we re-formulate the problem in rates and
eliminate thus the plastic distortion, casting thus a %so-called
hypoplastic model. Such a rate formulation based on velocity and
elastic strain together with the plastic distor{{t}}ion rate
is sometimes used in engineering, although without any rigorous analysis.
A conceptual benefit from avoiding the plastic distortion is elimination
of discussions about an intermediate stress free configuration
arising from it, which is felt as a fictitious and physically meaningless.
The eliminated plastic strain can be ``reconstructed'' a-posteriori.
At this point, we involve higher gradients in the dissipative potential, which
allows for  a  rigorous mathematical analysis, together with
rigorous control of invertibility of elastic strain by the stored energy.

The mentioned energy  dissipation  balance is used eventually in
Section~\ref{sec-anal} to perform the analysis of the hypoplastic model by a
discretization in space (Galerkin method) together with some regularization.
In this way, existence of weak solutions is proved by a constructive
method, giving also some conceptual numerical algorithm.

To highlight the main attributes of the model and its mathematical
treatment, we present it as quasistatic with the acceleration (and related
inertial forces) neglected and (still nonconstant and evolving) mass density thus
occurring only in the bulk-load term. A lot of nontrivial analytical
technicalities, now well understood from compressible fluid
dynamics \cite{Feir04DVCF,Lion98MTFM}, are thus avoided. For the same reason,
we present the model isothermally. The second simplifying
assumption (but most frequently adopted in literature) is
nonpenetrability of the boundary (i.e.\ normal velocity zero), which 
allows also for considering fixed boundary even for  the  Eulerian
description. In applied sciences, a rough
approach to live with this nonpenetrable boundary is considering 
time-varying domains embedded into a fictitious fixed domain and let
the material being inhomogeneous, composed from the viscoelastoplastic solid
and a very soft one. In geodynamical modelling,
this trick is sometimes called the sticky-air approach.
Eventually, we will exploit suitable gradient theories to facilitate
the proof of existence of weak solutions to the hypoplastic model.
There seems to be a general agreement that large-strain models
ultimately needs some higher gradients to cope with geometrical
nonlinearitites. In engineering models, various gradients are
used to control internal length-scales. Two principle options
are usage of higher gradients in the conservative way (i.e.\ enhancing
the stored energy) or in the dissipative way (i.e.\ enhancing
the dissipation potential). In Section~\ref{sec-hypo}, we will accept the
latter option. Of course, various combinations of both options
can be considered, too. Sometimes, even a diffusion is added into the 
evolution rule of the deformation gradient \cite{BFLS18EWSE},
which seems only artificial if not used in a modified
form \cite[Remark~3]{Roub??VELS} where it might have an interpretation of
Brenner's stress diffusion \cite{Bren05KVT}; in creep (fluid) models
cf.\ e.g.\ \cite{BuFeMa19CCVR} or for an incompressible case also
\cite{BaBuMa21LDET,BMPS18PDEA,EiHoMi21LHSV,MPSS18TVRT},
although even this is considered disputable.

The main notation used in this paper is summarized in the following table:

{\small
\begin{center}
\hspace*{1em}\begin{minipage}[t]{39em}
\hrule
\ \end{minipage}
\end{center}\vspace*{-1.5em}
\begin{center}
\ \ \ \ \begin{minipage}[t]{16em}%{0.25\linewidth}

$\vv$ velocity (in m/s),

$\varrho$ mass density (in kg/m$^3$),

$\FF$ deformation gradient,

$\Fe$ elastic strain,

$\Fp$ inelastic (plastic) distortion,

$\TT$ Cauchy stress (symmetric - in Pa),

$\SS=\varphi(\Fe)$ Piola stress (in Pa),

%$-\FF^\top\SS$ Eshelby stress (in Pa),

$\mathfrak{H}$ hyperstress (in Pa\,m),

\vspace*{.3em}

$\R_{\rm dev}^{d\times d}=\{A\in\R^{d\times d};\ {\rm tr}\,A=0\}$,

\end{minipage}\hspace*{1em}
\begin{minipage}[t]{22em}%{0.55\linewidth}

$\varphi=\varphi(\Fe)$ stored energy (in J/m$^3$=Pa),

$\ee(\vv)=\frac12\Nabla\vv^\top\!+\frac12\Nabla\vv$ small strain rate (in s$^{-1}$),

$\xi=\xi(\ee(\vv))$ viscous dissipation potential,

$\zeta=\zeta(\Lp)$ plastic dissipation potential,

$\Lp=\DTFp\Fpp^{-1}$ plastic distortion rate (in s$^{-1}$),

$(^{_{_{\bullet}}})\!\DT{^{}}=\pdt{}{^{_{_{\bullet}}}}+(\vv{\cdot}\Nabla)^{_{_{\bullet}}}$
convective time derivative,

$\GRAVITY$ external bulk load ({gravity}  acceleration in m/s$^{2}$),

$\ff$ traction load,

$\R_{\rm sym}^{d\times d}=\{A\in\R^{d\times d};\ A^\top=A\}$.

\end{minipage}
\end{center}
\vspace*{-1.2em}\begin{center}
\hspace*{1em}\begin{minipage}[t]{39em}
\hrule
\ \end{minipage}
\end{center}
}

\section{Plasticity at large strains classically}\label{sec-classic}
%        ~~~~~~~~~~~~~~~~~~~~~~~~~~~~~~~~~~~~~~
In large-strain continuum mechanics, the basic geometrical concept is the
time-evolving deformation $\yy:\varOmega\to\R^d$ as a mapping from a reference
configuration $\varOmega\subset\R^d$ into a physical space $\R^d$.
 The ``Lagrangian'' space variable in the reference configuration will be
denoted as $\XX\in\varOmega$ while in the ``Eulerian'' physical-space
variable by $\xx\in\R^d$.  The basic geometrical object is the deformation
gradient $\FF=\Nabla_{\XX}^{}\yy$.

 We will be interested in deformations $\xx=\yy(t,\XX)$ evolving in time,
which are sometimes called ``motions''. The important quantity is the
Eulerian  velocity $\vv=\DT\yy=\pdt{}\yy+(\vv{\cdot}\Nabla)\yy$.
Here and thorough the whole article, we use the dot-notation
$(\cdot)\!\DT{^{}}$ for the {\it convective time derivative}
applied to scalars or, component-wise, to vectors or tensors.

Then the velocity gradient
$\Nabla\vv=\nabla_{\!\XX}^{}\vv\nabla_{\!\xx}^{}\XX=\DT\FF\FF^{-1}$,
where we used the chain-rule calculus  and
$\FF^{-1}=(\nabla_{\!\XX}^{}\xx)^{-1}=\nabla_{\!\xx}^{}\XX$. 
This gives the {\it transport equation-and-evolution  for the
deformation gradient} as
\begin{align}
\DT\FF=(\nabla\vv)\FF\,.
\label{ultimate}\end{align}
From this, we also obtain the transport equation for the determinant
$\det\FF$ as
\begin{align}
\DT{\overline{\det\FF}}=(\det\FF)({\rm div}\,\vv).
\label{ultimate-det}\end{align}
 The understanding of \eqref{ultimate} and \eqref{ultimate-det} is a bit
delicate because it mixes the Eulerian $\xx$ and the Lagrangian $\XX$; note that
$\nabla\vv=\nabla_\xx\vv(\xx)$ while standardly $\FF=\nabla_{\!\XX}^{}\yy=\FF(\XX)$.
In fact, we consider $\FF{\circ}\bm\xi$ where $\bm\xi:\xx\mapsto\yy^{-1}(t,\XX)$
is the so-called {\it return} (sometimes called also a {\it reference})
{\it mapping}. Thus $\FF$ depends on $\xx$ and \eqref{ultimate} and
\eqref{ultimate-det} are equalities for a.a.\ $\xx$. The reference mapping
$\bm\xi$, which is well defined
through its transport equation $\DT{\bm\xi}=\bm0$, actually does not explicitly
occur in the formulation of the problem. Here we will benefit from the boundary
condition $\vv{\cdot}\nn=0$ below, which causes that the actual domain
$\varOmega$ does not evolve in time. The same concerns $\bm{T}$ in
\eqref{Euler1=plast} below, which will make the problem indeed fully Eulerian, as
announced in the title itself. Cf.\ the continuum-mechanics textbooks as e.g.\
\cite{GuFrAn10MTC,Mart19PCM}.

The mass density (in kg/m$^3$) is an extensive variable, and its transport
(expressing that the conservation of mass) writes as the {\it continuity equation}
$\pdt{}\varrho+{\rm div}(\varrho\vv)=0$, or, equivalently, the {\it mass transport
equation}
\begin{align}\label{cont-eq+}
\DT\varrho=-\varrho\,{\rm div}\,\vv\,.
\end{align}

Introducing a (generally non-symmetric) {\it plastic %deformation
distor{{t}}ion} tensor $\Fp$,
a conventional large-strain plasticity is based on Kr\"oner-Lie-Liu
%Green-Naghdi \cite{GreNag65GTEP}
\cite{Kron60AKVE,LeeLiu67FSEP} {\it multiplicative
decomposition}\index{multiplicative decomposition}
\begin{align}\label{Green-Naghdi}
\FF=\Fe\Fp\,.
\end{align}
The interpretation of $\Fp$ is a transformation of the reference
configuration into an intermediate stress-free configuration, and
then the {\it elastic strain}\index{elastic strain} $\Fe$ transforms this
intermediate configuration into the current actual configuration.

The main ingredients of the model are the (volumetric) {\it stored 
energy} and the {\it dissipation potential}, i.e.\ the physical unit of the
stored energy is  Pa=J/m$^3$ and of the dissipation potential is Pa/s.
The stored energy $\widehat\varphi(\FF,\Fp)$ depends naturally on the elastic strain
$\Fe=\FF\Fpp^{-1}$ and possibly also on $\Fp$ itself  if an isotropical
hardening were considered, but not directly on $\FF$.
In this section we will consider $\widehat\varphi(\FF,\Fp)=
\varphi(\FF\Fpp^{-1})$.
The other ingredient is the dissipation potential depending  on the
symmetric velocity gradient $\ee(\vv)=\frac12\Nabla\vv^\top\!+\frac12\Nabla\vv$
and on the {\it plastic distortion rate} $\DTFp\Fpp^{-1}$. We will consider this
dissipation potential as $\xi(\ee(\vv))+\widehat\zeta(\Fp,\DTFp)$
 with the plastic dissipation potential $\widehat\zeta$ depending on the
plastic distortion rate, i.e.\ $\widehat\zeta(\Fp,\DTFp)=\zeta(\DT\FF\Fpp^{-1})$ for
some potential $\zeta$. If quadratic, these two parts of the dissipation
potential involve linear Kelvin-Voigt-type and Maxwell-type viscosities into
the model, and altogether with the elastic part determined by the stored
energy, we obtain the {\it Jeffreys viscoelastic rheological model} in the
shear while the volumetric response obeys the {\it Kelvin-Voigt
rheology} if $\Fp$ is purely isochoric, as in Sections~\ref{sec-hypo} and
\ref{sec-anal} below. A quadratic $\zeta(\cdot)$ thus describes {\it creep}.
Yet, $\zeta$ may be non-quadratic and even non-smooth at the rate zero, which models an
``activated creep'' as in ice or {\it plasticity}, or even out of zero rate as in
the Tresca plasticity. This nonsmoothness makes the convex subdifferential
$\pl\zeta$ set-valued and thus why we wrote an inclusion ``$\ni$'' in \eq{Euler3=plast}.

The quasistatic evolution system then consists from the mass transport equation,
momentum equilibrium, the deformation gradient transport \eq{ultimate}, and
a flow rule for the plastic distortion $\Fp$. Specifically, in terms of
$\widehat\varphi$ and $\widehat\zeta$ the  system for $(\varrho,\vv,\FF,\Fp)$ reads as:
\begin{subequations}\label{Euler-plast}\begin{align}\label{Euler0=plast}
&\DT\varrho=-\varrho{\rm div}\,\vv\,,
\\&    {\rm div}\TT
    +\varrho\GRAVITY={\bm0}\,\ \ \text{ with }\ \
    \TT=\widehat\varphi_\FF'(\FF,\Fp)\FF^\top 
    +\widehat\varphi(\FF,\Fp)\bbI+\xi'(\ee(\vv))
\label{Euler1=plast}
\\\label{Euler2=plast}
&\DT\FF=(\Nabla\vv)\FF\,,
\\\label{Euler3=plast}
&\pl_{\DTFp}\widehat\zeta(\Fp,\DTFp)\ni-\widehat\varphi_{\Fp}'\!(\FF,\Fp)%\Fpp^\top
\,,
\end{align}\end{subequations}
where $\widehat\varphi_\FF'(\FF,\Fp)$ is the so-called {\it Piola stress}  and while
\eq{Euler3=plast} has the standard structure of the so-called Biot equation.

In \eq{Euler1=plast}, $\GRAVITY$ means a given acceleration (typically the {\it
gravity acceleration}) while we neglected the inertial force $\varrho\DT\vv$.
This last point substantially simplifies the analytical arguments below while
keeping the main phenomena under our focus in the game, although
the absence of the kinetic energy makes estimation of the bulk force
quite technical, cf.\ \eq{est-of-rho.f.v} below. In particular,
although \eq{Euler1=plast} neglects the acceleration $\DT\vv$ and thus the
mentioned inertial force $\varrho\DT\vv$, the mass $\varrho$ and its
transport \eq{Euler0=plast} are still involved. Such models are called
{\it quasistatic} (or, in geophysics, sometimes also {\it quasidynamic}).

It should be noted that the system \eq{Euler-plast} is truly standard,
and can often be found in literature, at least in its parts. Its
structure is, to a large extent, dictated by pursuing a consistent energetics
and the gradient doubly-nonlinear structure. The evolution of $\Fe$ and the
multiplicative decomposition \eq{Euler2=plast} is indeed most often
considered  as a  model for large-strain elastoplasticity and does
not need any comments here, as well as the mass transport \eq{Euler0=plast}
and the momentum equilibrium with the Kelvin-Voigt type Cauchy stress.
The  conservative ({elastic}) part of the Cauchy stress
$\widehat\varphi_\FF'(\FF,\Fp)\FF^{\top}+\widehat\varphi(\FF,\Fp)\bbI$ involves also a
pressure contribution since the free energy $\varphi$ is here considered
per  actual  volume (and not per  the referential volume or
 mass), cf.\ \cite[Rem.\,2]{Roub??VELS}; the notation $\bbI$ here
and in what follows stands for the identity matrix. 
The symmetry of such Cauchy stress is a standard consequence of the frame
indifference of $\varphi$ which is to be assumed, although we will not
explicitly use it. The form of the rate $%\Lp=
\DTFp\Fpp^{-1}$ which occurs in
the dissipation potential, is used most often in a position of an inelastic
distor{{t}}ion rate, cf.\
\cite{BesGie94MMID,ClTMau00ESTF,Dafa84PSCS,GurAna05DMSP,GuFrAn10MTC,Lee69EPDF,MiRoSa18GERV,RajSri04TMMN,XiBrMe00CFEP,ZNSM21MFSC},
and is also compatible with the so-called plastic indifference, cf.\ e.g.\
\cite{Miel03EFME}. The plastic flow rule \eq{Euler3=plast} is exactly as in
\cite{RouSte19FTCS}, cf.\ also \cite[Sec.9.4]{KruRou19MMCM}. Sometimes, however,
the plastic flow rule is formulated in the rate $\Fe\DTFp\Fpp^{-1}\Fee^{-1}$,
cf.\ Remark~\ref{rem-alternative} below or, in Lagrangian setting, as
$(\Fpp^\top\Fp)\!\DT{^{}}=\DTFp^\top\Fp+\Fpp^\top\DTFp$ in \cite{GraSte17FPPP},
too.

To reveal the mentioned energy  dissipation  balance behind the system
\eq{Euler-plast}, we
complete it by suitable boundary conditions. An important aspect is to impose
impenetrability of the boundary, which allows also for working on
a fix domain even in the Eulerian description. This also simplifies many
analytical arguments and is most often used in literature, too.
Moreover, in our quasistatic case where inertial forces are neglected, we need
to fix the body at a part of the boundary at least viscously. Thus we
consider a combination of a homogeneous Dirichlet combination
in the normal direction and the Newton (or Navier) condition
in the tangential direction:
\begin{align}\label{Euler-plast-BC}
\vv{\cdot}\nn=0\ \ \ \text{ and }\ \ \ [\TT\nn]_\text{\sc t}^{}
+\kappa\vv_\text{\sc t}^{}=\ff\,,
\end{align}
where $(\cdot)_\text{\sc t}^{}$ denotes the tangential component of a vector
on the boundary $\varGamma$ and $\nn$ is the unit outward normal to $\varGamma$.
The first condition in \eq{Euler-plast-BC} simplifies considerably
the situation and allows for working on a fixed domain $\varOmega$. Then, formally,
we obtain the energetics by testing \eq{Euler1=plast} by $\vv$ and using
\eq{Euler2=plast} tested by %$\SS\Fpp^{-\top}$
$\SS$ and  by testing \eq{Euler3=plast} by
%$\Lp$
$\DTFp\Fp^{-1}$, while \eq{Euler0=plast} does not directly contribute to
the energetics because the inertial term $\varrho\DT\vv$ has been
neglected. The former test gives
\begin{align}\nonumber
&
\int_\varOmega{\rm div}\,\TT{\cdot}\vv\,\d x=\int_\varGamma(\TT\nn){\cdot}\vv\,\d S
-\int_\varOmega\TT{:}\ee(\vv)\,\d x
\\[-.1em]&\nonumber=\!
\int_\varGamma
(\TT\nn){\cdot}\vv\,\d S-\!\int_\varOmega
\Big(\widehat\varphi_{\FF}'(\FF,\Fp)\FF^\top\!
+\widehat\varphi(\FF,\Fp)\bbI+\xi'(\ee(\vv))%+{\bm S}_{\rm str}
\Big){:}\ee(\vv)\,\d x
\\[-.1em]&\nonumber=\!
\int_\varGamma(\TT\nn){\cdot}\vv\,\d S-\!\int_\varOmega
\Big(\widehat\varphi_{\FF}'(\FF,\Fp){:}(\nabla\vv)\FF
+\widehat\varphi(\FF,\Fp){\rm div}\,\vv+\xi'(\ee(\vv)){:}\ee(\vv)\,\d x
\\[-.1em]&\nonumber=\!\int_\varGamma\!(\TT\nn){\cdot}\vv\,\,\d S
-\!\int_\varOmega\!\widehat\varphi_\FF'(\FF,\Fp){:}\DT\FF
+\widehat\varphi(\FF,\Fp){\rm div}\,\vv
+\xi'(\ee(\vv)){:}\ee(\vv)\,\d x
\\[-.1em]&\nonumber=
\!\int_\varGamma\!(\TT\nn){\cdot}\vv\,\,\d S
-\frac{\d}{\d t}\int_\varOmega\widehat\varphi(\FF,\Fp)\,\d x
\\[-.4em]&\nonumber\hspace{2em}
-\int_\varOmega\widehat\varphi_\FF'(\FF,\Fp){:}(\vv{\cdot}\Nabla)\FF
+\widehat\varphi(\FF,\Fp){\rm div}\,\vv
+\xi'(\ee(\vv)){:}\ee(\vv)
-\widehat\varphi_{\Fp}'(\FF,\Fp){:}\pdt{\Fp\!}\,\d x
\\[-.1em]&\nonumber=
\!\int_\varGamma\!(\TT\nn){\cdot}\vv\,\,\d S
-\frac{\d}{\d t}\int_\varOmega\widehat\varphi(\FF,\Fp)\,\d x
-\int_\varOmega\Big(\widehat\varphi_\FF'(\FF,\Fp){:}(\vv{\cdot}\Nabla)\FF
+\widehat\varphi(\FF,\Fp){\rm div}\,\vv
\\[-.2em]&\nonumber\hspace{11em}
+\xi'(\ee(\vv)){:}\ee(\vv)
-\widehat\varphi_{\Fp}'(\FF,\Fp){:}\DTFp
+\widehat\varphi_{\Fp}'(\FF,\Fp){:}(\vv{\cdot}\Nabla)\Fp\Big)\,\d x
\\[-.1em]&=
\!\int_\varGamma\!\ff{\cdot}\vv{-}\kappa|\vv|^2\,\d S
-\frac{\d}{\d t}\!\int_\varOmega\widehat\varphi(\FF,\Fp)\,\d x
-\!\int_\varOmega\!\xi'(\ee(\vv)){:}\ee(\vv)
+\pl_{\DTFp}\!\widehat\zeta(\Fp,\DTFp){:}\DTFp\,\d x\,,
\label{Euler-large-1}\end{align}
where the last equality results when using the inclusion \eq{Euler3=plast} tested by
% the plastic distortion rate $\DTFp\Fpp^{-1}$.
$\DTFp$; actually, such test generally gives only an inequality but we implicitly
rely on certain smoothness of $\zeta$ out of $0$ as assumed later in \eq{Euler-ass-zeta}.
Here we have used several times the matrix algebra 
\begin{align}\label{algebra}
A:(BC)=(B^\top A):C=(AC^\top):B
\end{align}
for any three square matrices $A$, $B$, and $C$. For the last equality in
\eq{Euler-large-1}, we also used the calculus
\begin{align}\nonumber
&\int_\varOmega\widehat\varphi_\FF'(\FF,\Fp){:}(\vv{\cdot}\Nabla)\FF
+\widehat\varphi_{\Fp}'(\FF,\Fp){:}(\vv{\cdot}\Nabla)\Fp
+\widehat\varphi(\FF,\Fp){\rm div}\,\vv\,\d x
\\&\qquad
=\int_\varOmega\nabla\widehat\varphi(\FF,\Fp){\cdot}\vv
+\widehat\varphi(\FF,\Fp){\rm div}\,\vv\,\d x
=\!\int_\varGamma\widehat\varphi(\FF,\Fp)(\vv{\cdot}\nn)\,\d S=0\,,
\label{Euler-pressure-calculus}\end{align}
where we employed the Green formula and the boundary conditions $\vv{\cdot}\nn=0$.
Altogether, we obtain (at least formally) the expected energy  dissipation  balance
\begin{align}\nonumber
  &\!\!\!\!\!\frac{\d}{\d t}
  \int_\varOmega\!\!\!\!\linesunder{\widehat\varphi(\FF,\Fp)%{+}\phi(\Fp)
  \!_{_{_{_{_{_{}}}}}}}{stored}{energy}\!\!\!\!\d x
+\!\int_\varOmega\!\!\!\!\lineunder{\xi'(\ee(\vv)){:}\ee(\vv)
+\pl_{\DTFp}\widehat\zeta(\Fp,\DTFp){:}\DTFp%\Fpp^{-1}
\!}{bulk dissipation rate}\!\!\!\d x
\\[-.6em]&\hspace{9em}
+\!\int_\varGamma\!\!\!\!\!\!\!\!\!\!\linesunder{\kappa|\vv|^2}{boundary}{dissipation rate}\!\!\!\!\!\!\!\!\!\!\d S
=\int_\varOmega\!\!\!\!\!\!\!\!\!\linesunder{\varrho\,\GRAVITY{\cdot}\vv_{_{_{_{_{}}}}}\!\!\!}{power of}{external load}\!\!\!\!\!\!\!\!\d x+\!\int_\varGamma\!\!\!\!\!\!\!\!\!\!\!\linesunder{\ff{\cdot}\vv_{_{_{_{_{}}}}}}{power of}{traction load}\!\!\!\!\!\!\!\!\!\!\d S\,.
  \label{energy++++}\end{align}

Actually, a usual assumption is that
the inelastic deformation $\Fp$ only concerns shear and does not
affect volume variations. We call it {\it isochoric} and it means
not only that that $\det\Fp$ is positive (to make $\Fp$ invertible) but even
that $\det\Fp=1$. Yet, the constraint $\det\Fp=1$ is not affine and, if
it would be ensured by the conservative part (i.e.\ through the stored energy
$\phi$) and thus explicitly involved into \eq{Euler-plast} together with a
corresponding Lagrange multiplier, the analytical treatment of such a
differential-algebraic-type system would become extremely difficult and likely
impossible.  One modelling option is to consider this isochoric constraint only
approximately by casting a hardening-like term acting on $\det\Fp$ to ensure that
$\det\Fp$ is positive and close to 1;
cf.\ \cite{DaRoSt??NHFV,KruRou19MMCM,MieRou16RIEF,RouSte19FTCS} for a Lagrangian
formulation.  Another option is, instead of the control of $\det\Fp$ in the stored
energy (implementing thus the isochoricity only approximately),
to implement the isochoricity exactly in the dissipative
part  relying on the calculus
\begin{align}\label{D/DT-of-determinant}
\DT{\overline{\det\Fp}}={\rm Cof}\Fp{:}\DTFp=(\det\Fp)\Fp^{-\top}{:}\DTFp=(\det\Fp){\rm tr}(\DTFp\Fpp^{-1})
\end{align}
and  by considering
\begin{align}\label{ass-inelastic-dissip}
\zeta:\R^{d\times d}\to[0,+\infty]\ \text{ is convex},\ \ \zeta(0)=0,\ \
\zeta\big(\R^{d\times d}{\setminus}\R_{\rm dev}^{d\times d}\big)=+\infty
\end{align}
together with precribing the contraint $\det\Fp=1$ on the initial condition;
cf.\ also  \cite[Sect.\,91.3]{GuFrAn10MTC}.  
Then \eq{Euler3=plast} ensures ${\rm tr}(\DTFp\Fpp^{-1})=0$ and, 
by \eq{D/DT-of-determinant}, $\det\Fp=1$
provided the initial inelastic deformation is isochoric. It is important
that the trace-free constraint %${\rm tr}\Lp=0$
in \eq{ass-inelastic-dissip} is linear, in contrast to
the non-affine constraint $\det\Fp=1$. 

Let us however emphasize that the rigorous analysis of the system
\eq{Euler-plast} would need still gradients in the dissipation potential,
which we will use in the following sections, cf.\ \eq{dissip-pot-grad} below,
but which we intentionally ignored in \eq{Euler-plast} in order to explain
the main structure of the model without unnecessary technicalities.

\begin{remark}[{\sl A gradient structure of $(\vv,\FF)$}]\upshape
 Implicitly, \ we have in mind the situation when $(\vv,\FF)$ is a
gradient in the sense that $(\vv,\FF)=[(\cdot)\!\DT{^{\,}},\!\Nabla]\yy$ of some deformation
$\yy$ which, however, does not explicitly occur in \eq{Euler-plast}. Indeed, an
existence of some
$\yy$ so that $\FF=\nabla_{\!\XX}^{}\yy$ and $\vv=\DT\yy$ is not automatic even
if $\FF$ is a gradient of some deformation at an initial time. Rather,
we can always construct the return
mapping $\bm\xi$ mentioned above by solving the simple transport equation
$\DT{\bm\xi}=\bm0$ with the initial condition $\bm\xi(0)$=identity. Then
$\FF=(\nabla_\xx\bm\xi)^{-1}$ and, if $\bm\xi(t):\varOmega\to\varOmega$ is
injective, the underlying deformation is $\yy(t)=\bm\xi^{-1}(t)$. This global
injectivity seems not automatic, however; cf.\ also \cite[Rem.\,7]{Roub??VELS}.
\end{remark}

\section{Rate form of plasticity -- hypoplasticity}\label{sec-hypo}
%        ~~~~~~~~~~~~~~~~~~~~~~~~~~~~~~~~~~~~~~~~~~~~~~~

We will now express the original model \eq{Euler-plast}  in terms
of the energies $\varphi$ and $\zeta$ instead of $\widehat\varphi$ and
$\widehat\zeta$. By this  way, % that
the plastic distortion $\Fp$ will be eliminated from the model, although
it will be possible to reconstruct it if the initial condition
is known. The plasticity evolution will be formulated exclusively
in terms of plastic distortion rate $\Lp$, cf.\ \eq{dafa-formula} below, 
called  here  {\it hypoplasticity}  in parallel how a rate
formulated hyperelasticity is called hypo-elasticity \cite{Treu55HE}. 
Actually, formulating the model in terms of $\Fe$ and $\Lp$ instead of the
multiplicative decomposition, was explicitly advocated
in \cite[p.249]{BesGie94MMID}, emphasizing that $\Fp$ bears no physical
relevance.

Moreover, we will %be able also to
eliminate the equation for the mass transport \eq{Euler0=plast},
although it will stay implicitly contained in the model. Relying on \eq{ultimate-det},
one can determine the density $\varrho$ instead of the transport equation for mass
density \eq{cont-eq+} from the algebraic relation
\begin{align}
\varrho=\frac{\varrho_0}{\det\FF}
\label{density-algebraically}\end{align}
where $\varrho_0$ is the mass density in the  reference configuration.
Indeed, one has the calculus
\begin{align}
\frac{\DT\varrho}{\varrho}=
\Bigg(\varrho_0\,\,\DT{\overline{\!\!\bigg(\frac1{\det\FF}\bigg)\!\!}}\,\,
+\frac{\DT\varrho_0}{\det\FF}\Bigg)
\frac{\det\FF}{\varrho_0}=-\frac{\DT{\overline{\det\FF}}}{\det\FF}
=-{\rm div}\,\vv
\label{towards-cont-eq}\end{align}
because $\DT\varrho_0=0$ and because,
 analogously to \eq{D/DT-of-determinant},  we have
\begin{align}\nonumber
\label{DT-det}\DT{\overline{\det\FF}}&=
{\rm Cof}\FF{:}\DT\FF=(\det\FF)\FF^{-\top}\!{:}\DT\FF
\\[-.3em]&=(\det\FF)
\bbI{:}\DT\FF\FF^{-1}=(\det\FF)\bbI{:}\Nabla\vv=(\det\FF){\rm div}\,\vv\,;
\end{align}
here we used also \eq{ultimate} and the matrix algebra \eq{algebra}.
Thus, the last equality in \eq{towards-cont-eq} is the transport equation
\eq{ultimate-det} while \eq{towards-cont-eq}
itself is just the continuity equation %\eq{Euler0+++}
\eq{cont-eq+}.  This would allow (and is actually often used) for
elimination of the continuity equation \eq{Euler0=plast} in
Section~\ref{sec-classic}.

 Here, assuming again \eq{ass-inelastic-dissip} and isochoricity of the
initial plastic distortion and,  thus, having $\det\Fp=1$ during the whole
evolution, we have $\det\FF=\det(\Fe\Fp)=\det\Fe\det\Fp=\det\Fe$ and 
\eq{density-algebraically} can be written as
\begin{align}
\varrho=\frac{\varrho_0}{\det\Fe}\,.
\label{density-algebraically+}\end{align}

%Third
 Moreover, applying the material derivative on
\eq{Green-Naghdi} and using \eq{Euler2=plast}, we obtain
$(\Nabla\vv)\FF=\DT\FF=\DTFe\Fp+\Fe\DTFp$ and, multiplying it by
$\FF^{-1}=\Fpp^{-1}\Fee^{-1}$, eventually we obtain
\begin{align}
\Nabla\vv=\!\!\!\!\!\!\morelinesunder{\DTFe\Fee^{-1}\!\!\!}{elastic\ \ }{distortion}{rate}\!\!\!+\Fe\!\!\!\!\!\!\!\morelinesunder{\DTFp\Fpp^{-1}\!\!\!\!}{plastic}{distortion}{rate $=:\Lp$}\!\!\!\!\Fee^{-1},
\label{dafa-formula}\end{align}
cf.\ e.g.\ \cite{BesGie94MMID,BeReBo16FSEE,Dafa84PSCS,GurAna05DMSP,GuFrAn10MTC,Hash20NCMF,HasYam13IFST,JirBaz02IAS,KhaHua95CTP,Lee69EPDF,MauEps98GMSE,RajSri04TRC,RajSri04TMMN,XiBrMe00CFEP,ZNSM21MFSC};
the terms ``distortion rates'' are due to \cite{GurAna05DMSP,GuFrAn10MTC}
while sometimes $\Lp$ is called a ``plastic dissipation tensor''
\cite{BesGie94MMID} or  ``velocity gradient of purely plastic deformation''
in \cite{Lee69EPDF}, etc.

Interestingly, in terms of plastic distortion rate $\Lp$, we will not see
explicitly the plastic distortion $\Fp$ and, multiplying \eq{dafa-formula} by
$\Fe$, we obtain an evolution rule for $\Fe$ even without any explicit
occurrence of $\Fee^{-1}$, namely
\begin{align}
\DTFe=(\Nabla\vv)\Fe-\Fe\Lp\,.
\label{evol-of-E}\end{align}

Mainly for analytical reasons, we will enhance the dissipation potential from
Sect.~\ref{sec-classic} by generally non-quadratic gradient terms as
\begin{align}\label{dissip-pot-grad}
  \xi(\ee(\vv))+\zeta(\Lp)+\frac{\NU}p|\Nabla\ee(\vv)|^p+\frac{\MU}q|\Nabla\Lp|^q
\end{align}
with some  (presumably small) coefficients  $\NU,\MU>0$,
cf.\ Remark~\ref{rem-grad} below. In the next Section~\ref{sec-anal}, we will
need both the gradient-term exponents sufficiently big, namely $p>d$ and $q>d$. 

 The stress $\widehat\varphi_\FF'(\FF,\Fp)\Fp^\top$ in \eq{Euler1=plast}
is to be written in terms of $\varphi$ by the calculus
$$
\widehat\varphi_\FF'(\FF,\Fp)\FF^\top=
\big[\varphi(\FF\Fp^{-1})\big]_\FF'\FF^\top=
\varphi'(\Fe)\Fp^{-\top}(\Fe\Fp)^\top=\varphi'(\Fe)\Fe^\top.
$$ 
When writ{ing}  the plastic flow rule \eq{Euler3=plast} in terms of
$\Lp$ as a purely ``algebraic'' relation without  any  explicit
reference to  the relation  $\Lp=\DTFp\Fpp^{-1}$ and using
\eq{density-algebraically} together with the
isochoric-inelasticity concept, we obtain the hypo-elastoplastic system as
\begin{subequations}\label{Euler-hypoplast}
\begin{align}
\nonumber
&{\rm div}{\bm\varSigma}+\varrho\GRAVITY={\bm0}\ \ \ \,\text{ with }\ \ \varrho=\frac{\varrho_0}{\det\Fe\!\!}\ ,
     \ \ \ {\bm\varSigma}=\TT-{\rm div}\,\mathfrak{H},
\
\\\nonumber
    &\hspace*{8em}\text{and }\ \ \ %{\bm\varSigma}
\TT=\SS\Fee^\top\!+\varphi(\Fe)\bbI+\xi'(\ee(\vv))
    \\
    &\hspace*{8em}\text{where }
    \SS=\varphi'(\Fe)  \ \ \text{ and }\ \ 
\mathfrak{H}=\NU|\nabla\ee(\vv)|^{p-2}\nabla\ee(\vv)\,,
\label{Euler1=hypoplast}
\\\label{Euler2=hypoplast}
&\DTFe=(\Nabla\vv)\Fe-\Fe\Lp\,,
\\\label{Euler3=hypoplast}
&\pl\zeta(\Lp)-{\rm div}(\MU|\Nabla\Lp|^{q-2}\Nabla\Lp)\ni\Fee^{\top}\SS\,.
\end{align}\end{subequations}
The right-hand side in \eq{Euler3=hypoplast}, being a driving stress
for the plastification process, can be identified as the
{\it Eshelby stress} \cite{BeReBo16FSEE,ClTMau00ESTF,EpsMau10RUES,MauEps98GMSE};
actually, the Eshelby stress standardly contains also a pressure part like 
does also  the elastic Cauchy stress
$\varphi'(\Fe)\Fe^{\top}\!+\varphi(\Fe)\bbI$ but such a pressure would not
affect the isochoric plastic evolution.  Again, the form of this stress is
dictated essentially in order to achieve the desired energy  dissipation
 balance \eq{energy++++}, i.e.\ now \eq{energy-hypoplast} below.

The system \eq{Euler-hypoplast} is to be completed by suitable
boundary conditions counting also the gradient terms arising from
the enhanced dissipation potential \eq{dissip-pot-grad}, say
\begin{align}\label{Euler-hypoplast-BC}
&\vv{\cdot}\nn=0\,,\ \ 
\big[{\bm\varSigma}\nn{-}\divS(\mathfrak{H}\nn)\big]_\text{\sc t}^{}\!
+\kappa\vv_\text{\sc t}^{}=\ff\,,
\ \ 
\Nabla\ee(\vv){:}(\nn{\otimes}\nn)={\bm0}\,,\ \text{ and }\ 
\Nabla\Lp{\cdot}\nn={\bm0}\,,%\ \ \text{ on }\ \varGamma,
\end{align}%\end{subequations}
where the $(d{-}1)$-dimensional surface divergence is defined as
\begin{align}\label{def-divS}
\divS={\rm tr}(\nablaS)\ \ \ \text{ with }\ \
\nablaS v=\nabla v-\frac{\partial v}{\partial\nn}\nn\,,
\end{align}
where ${\rm tr}(\cdot)$ is the trace of a $(d{-}1){\times}(d{-}1)$-matrix and
$\nablaS v$ is the surface gradient of $v$.

The energetics behind the model \eq{Euler-hypoplast} can be revealed 
by testing %\eqref{Euler0=hypoplast} by $\frac12|\vv|^2$,
\eqref{Euler1=hypoplast} by $\vv$, and \eqref{Euler2=hypoplast} by $\SS$ (or, more
precisely, testing \eq{dafa-formula} by $\SS\Fee^\top$), and
\eqref{Euler3=hypoplast} by $\Lp$. Using
the Green formula and %the matrix algebra \eq{algebra} together with
\eq{dafa-formula} tested by $\SS\Fee^\top$,
we obtain from the Cauchy stress $\TT$:
\begin{align}\nonumber
&\!\!\int_\varOmega\!\!{\rm div}\,\TT{\cdot}\vv\,\d x
=\!\int_\varGamma\!\vv{\cdot}\TT\nn\,\d S
-\!\!\int_\varOmega\!(\SS\Fee^\top){:}\Nabla\vv+\varphi(\Fe){\rm div}\,\vv
+\xi'(\ee(\vv)){:}\ee(\vv)\,\d x
\\&\nonumber
=\!\int_\varGamma\!\vv{\cdot}\TT\nn\,\d S
-\!\!\int_\varOmega\!
\varphi'(\Fe)\Fee^\top{:}\big(\DTFe\!{+}\Fe\Lp\big)\Fee^{-1}
+\varphi(\Fe){\rm div}\,\vv+\xi'(\ee(\vv)){:}\ee(\vv)\,\d x
\\&\nonumber
=\!\int_\varGamma\!\vv{\cdot}\TT\nn\,\d S
-\!\!\int_\varOmega\!
\varphi'(\Fe){:}\DTFe
+\Fee^\top\varphi'(\Fe){:}\Lp+\varphi(\Fe){\rm div}\,\vv
+\xi'(\ee(\vv)){:}\ee(\vv)\,\d x
\\&
=\!\int_\varGamma\!\vv{\cdot}\TT\nn\,\d S
-\frac{\d}{\d t}\!\int_\varOmega\varphi(\Fe)\,\d t
-\!\int_\varOmega\pl\zeta(\Lp){:}\Lp+\MU|\Nabla\Lp|^q
+\xi'(\ee(\vv)){:}\ee(\vv)\,\d x\,,
\label{Euler-hypoplast-test-momentum}\end{align}
where we also used the matrix algebra \eq{algebra} for 
$\varphi'(\Fe)\Fee^\top{:}(\DTFe\Fee^{-1})=
\varphi'(\Fe)\Fee^\top\Fee^{-\top}{:}\DTFe=\varphi'(\Fe){:}\DTFe$
and for $\varphi'(\Fe)\Fee^\top{:}(\Fe\Lp\Fee^{-1})=
\varphi'(\Fe){:}(\Fe\Lp)=\Fee^\top\varphi'(\Fe){:}\Lp$.
In comparison with Sect.\,\eq{sec-classic}, note that
$\pl_{\DTFp}\widehat\zeta(\Fp,\DTFp){:}\DTFp=\pl\zeta(\Lp){:}\Lp$. The pressure
term is to be treated similarly as in \eq{Euler-pressure-calculus} by the calculus
\begin{align}\nonumber
&\!\int_\varOmega\varphi'(\Fe)\bigg(\pdt{\Fe\!}
+(\vv{\cdot}\Nabla)\Fe\bigg)
+\varphi(\Fe)\,{\rm div}\,\vv\,\d x
\\&=\int_\varOmega\pdt{}\varphi(\Fe)+\nabla\varphi(\Fe){\cdot}\vv
+\varphi(\Fe){\rm div}\,\vv\,\d x
=\frac{\d}{\d t}\int_\varOmega\!\varphi(\Fe)\,\d x
+\!\int_\varGamma\varphi(\Fe)(\!\!\!\!\!\lineunder{\vv{\cdot}\nn}{$=0$}\!\!\!\!\!)\,\d S.\!
\label{Euler-hypoplast-test-momentum+}\end{align}
The further contribution from the hyperstress gives, using Green formula over
$\varOmega$ twice and the surface Green formula over $\varGamma$, that
\begin{align}\nonumber
\int_\varOmega{\rm div}^2\mathfrak{H}{\cdot}\vv\,\d x&=\int_\varGamma\vv{\cdot}
{\rm div}\mathfrak{H}\nn\,\d S-\int_\varOmega{\rm div}\mathfrak{H}{:}\Nabla\vv\,\d x
\\&\nonumber=\int_\varOmega\mathfrak{H}\Vdots\Nabla^2\vv\,\d x+\int_\varGamma
\nn{\cdot}\mathfrak{H}{:}\Nabla\vv-\vv{\cdot}
{\rm div}\mathfrak{H}\nn\,\d S
\\&\nonumber=\int_\varOmega\mathfrak{H}\Vdots\Nabla^2\vv\,\d x+\int_\varGamma
\mathfrak{H}{:}(\nn{\otimes}\nn)
+\nn{\cdot}\mathfrak{H}{:}\NablaS\vv-\vv{\cdot}
{\rm div}\mathfrak{H}\nn\,\d S
\\&=\int_\varOmega\NU|\nabla \ee(\vv)|^p\,\d x+\int_\varGamma
\mathfrak{H}{:}(\nn{\otimes}\nn)-
\big(\divS(\nn{\cdot}\mathfrak{H})+{\rm div}\mathfrak{H}\nn\big)
{\cdot}\vv\,\d S\,,
\label{Euler-hypoplast-test-momentum++}\end{align}
where we used the decomposition of $\Nabla\vv$ into its normal and tangential
parts, i.e.\ written componentwise $\nabla\vv_i=(\nn{\cdot}\nabla\vv_i)\nn+\nablaS\vv_i$.

Merging the boundary integrals in \eq{Euler-hypoplast-test-momentum}
and in \eq{Euler-hypoplast-test-momentum++} and using the boundary condition
$({\bm\varSigma}\nn-\divS(\mathfrak{H}\nn))_\text{\sc t}^{}
+\kappa\vv_\text{\sc t}^{}=\ff$,
we thus obtain
(at least formally) the energy   dissipation  balance
\begin{align}\nonumber
  &\frac{\d}{\d t}
  \int_\varOmega\!\!\!\!\linesunder{\varphi(\Fe)
  \!_{_{_{_{_{_{}}}}}}}{stored}{energy}\!\!\!\!\,\d x
%\\[-.2em]&\nonumber\hspace{3.5em}
  +\int_\varOmega\!\!\!\!\lineunder{\xi'(\ee(\vv)){:}\ee(\vv)
%+\GM|\Lp|^2
+\pl\zeta(\Lp){:}\Lp+\NU|\nabla \ee(\vv)|^p+\MU|\Nabla\Lp|^q_{_{_{}}}
\!}{dissipation rate}\!\d x
\\[-.9em]&\hspace{15em}
+\!\int_\varGamma\!\!\!\!\!\!\!\!\!\!\linesunder{\kappa|\vv|^2}{boundary}{dissipation rate}\!\!\!\!\!\!\!\!\!\!\d S
=\int_\varOmega\!\!\!\!\!\!\!\!\!\linesunder{\varrho\GRAVITY{\cdot}\vv_{_{_{_{_{}}}}}\!\!\!}{power of}{external load}\!\!\!\!\!\!\!\d x+\int_\varGamma\!\!\!\!\!\!\!\!\!\!\linesunder{\ff{\cdot}\vv_{_{_{_{_{}}}}}\!\!}{power of}{traction load}\!\!\!\!\!\!\!\!\!\!\d S\,.
  \label{energy-hypoplast}\end{align}
If one is interested in an ``a posteriori'' reconstruction of the plastic
distortion $\Fp$,  one should prescribe also an initial condition
$\Fp|_{t=0}^{}=\Fpzero$ and, by re-arranging \eq{dafa-formula}, to use the
{\it plastic-strain evolution} rule
\begin{align}
  \DTFp=\Lp\Fp\,.
\label{plastic-strain-reconstructed}\end{align}
Only at this ``a posteriori'' point, one should
consider the assumption $\det\Fpzero=1$ on which the system
\eq{Euler-hypoplast} relied when arising from \eq{Euler-plast}.

\begin{remark}[{\sl Gradient theories in rates}]\label{rem-grad}\upshape
So-called gradient theories in continuum mechanical models are
very standard, determining some internal length scales and often
facilitating mathematical analysis. They can be applied to the stored
energy or to the dissipative potential, i.e.\ they contribute to
the conservative or to the dissipative parts of the model.
Here we used the latter option in \eq{dissip-pot-grad}.
The first gradient term leads to the hyperstress $\mathfrak{H}$ in the
momentum equation \eq{Euler1=hypoplast} and is compliant with the 
so-called {\it 2nd-grade nonsimple fluid} concept devised
by E.\,Fried and M.\,Gurtin \cite{FriGur06TBBC}
  %\cite{FriGur05SGFT,FriGur06TBBC}
and, earlier and even more generally and nonlinearly, as {\it multipolar fluids} by
J.\,Ne\v cas at al.\
\cite{BeBlNe92PBMV,BeNeRa99EUFM,Neca94TMF,NeNoSi89GSIC,NecRuz92GSIV,NecSil91MVF},
following ideas of R.A.\,Toupin \cite{Toup62EMCS} and R.D.~Mindlin \cite{Mind64MSLE}
for elastic solids. The further gradient term in \eq{dissip-pot-grad} gives rise to
${\rm div}(\MU|\Nabla\Lp|^{q-2}\Nabla\Lp)$ in the plastic-rate
evolution \eq{Euler3=hypoplast}. This causes a certain ``dynamical'' hardening
involves a certain length scale to the plastic distortion but,
does not make spurious hardening effects during long lasting plastification
or creep, unlike the conventional kinematic or isotropic hardening in the
conservative part. Similarly, in Lagrangian formulation,
\cite{DaRoSt??NHFV} used the plastic distortion rate in $\Delta\pdt{}\Fp$.
\end{remark}

\begin{remark}[{\sl An alternative model}]\label{rem-alternative}\upshape
  There is not a general agreement on an interpretation of the
  additive  split of the rate in  \eq{dafa-formula}.
 One can also work with the plastic rate as
$\Lp=\Fe\DTFp\Fpp^{-1}\Fee^{-1}$, cf.\ \cite[Sect.10.4]{Bert08EPLD},
\cite[Formula (7.1.4)]{BesGie94MMID},
\cite[Formulae (2.5)--(2.7)]{Krat73FSTE}, 
or \cite[Formulae (95)--(96)]{Volo13AELD}. This 
gives the elastic-strain evolution rule
$\DTFe=(\Nabla\vv-\Lp)\Fe$ instead of \eq{evol-of-E}, cf.\
\cite[Sect.7.1]{BesGie94MMID}, and the right-hand side of the
plastic flow rule \eq{Euler3=hypoplast} is 
$-\SS\Fee^\top$ instead of the Eshelby stress $-\Fee^\top\SS$. The dissipation
potential should act on this alternative $\Lp$.
The plastic distortion can be reconstructed, instead of
\eq{plastic-strain-reconstructed}, by 
$\DTFp=\Fee^{-1}\Lp\Fe\Fp$. This alternative model can capture the zero
plastic spin, provided the plastic spin is understood as the skew-symmetric
part of $\Lp$ because the driving (Cauchy) stress $\SS\Fee^\top$ is
symmetric. Sometimes, the plastic spin is however understood as the
skew-symmetric part of $\Lp$ when the splitting \eq{dafa-formula} is used,
cf.\ \cite[Sect.\,91.3]{GuFrAn10MTC}.
For the discussion according both variant see \cite[Sect.7.1]{BesGie94MMID}
or \cite[Sect.4]{NeNa83FDEP}.
Actually, for small elastic deformations where $\Fe\sim\bbI$,
both variants do not  differ much  from each other.
\end{remark}

\section{Analysis -- weak solutions of \eq{Euler-hypoplast}}\label{sec-anal}
%        ~~~~~~~~~~~~~~~~~~~~~~~~~~~~~~~~~~~~~~~~~~~~~~~~~~

We will provide a proof of existence and certain regularity of weak solutions.
It should be emphasized that, even with the nonlinear dissipative gradient
terms which have regularizing effects, it is quite nontrivial. %\COMMENT{DETAILS??}
 
We will use the standard notation concerning the Lebesgue and the Sobolev
spaces, namely $L^p(\varOmega;\R^n)$ for Lebesgue measurable functions
$\varOmega\to\R^n$ whose Euclidean norm is integrable with $p$-power, and
$W^{k,p}(\varOmega;\R^n)$ for functions from $L^p(\varOmega;\R^n)$ whose
all derivative up to the order $k$ have their Euclidean norm integrable with
$p$-power. We also write briefly $H^k=W^{k,2}$. The notation
$p^*$ will denote the exponent from the embedding
$W^{1,p}(\varOmega)\subset L^{p^*}(\varOmega)$, i.e.\ $p^*=dp/(d{-}p)$
for $p<d$ while $p^*\ge1$ arbitrary for $p=d$ or $p^*=+\infty$ for $p>d$.
Moreover, for a Banach space
$X$ and for $I=[0,T]$, we will use the notation $L^p(I;X)$ for the Bochner
space of Bochner measurable functions $I\to X$ whose norm is in $L^p(I)$
while $W^{1,p}(I;X)$ denotes for functions $I\to X$ whose distributional derivative
is in $L^p(I;X)$. Also, $C(\cdot)$ and $C^1(\cdot)$
will denote spaces of continuous and continuously differentiable functions.

A highly applicable assertion was originally
devised for situations when $\FF=\Nabla\yy$ with $\yy\in W^{2,p}(\varOmega;\R^d)$
but actually it holds in more general situations, as used also in
\cite{KruRou19MMCM,MieRou16RIEF,MieRou20TKVR,RouSte19FTCS}:

\begin{lemma}[T.J.\,Healey and S.\,Kr\"omer \cite{HeaKro09IWSS}]\label{lem-Healey-Kromer}
Let $\varkappa>rd/(r{-}d)$ for
some $r>d$. Then, for any $C<+\infty$, there is $\epsilon>0$ such
that, for any $\FF\in W^{1,r}(\varOmega;\R^{d\times d})$ with $\det\FF>0$
a.e.\ on $\varOmega$, it holds 
$$
\|\FF\|_{W^{1,r}(\varOmega;\R^{d\times d})}^{}+\int_\varOmega\frac1{(\det\FF)^\varkappa\!}\,\d x
\le C\qquad\Rightarrow\qquad\det\FF\ge\epsilon\ \text{ on }\ \barOmega\,.
$$
\end{lemma}

To devise a weak formulation of the initial-boundary-value problem for
the system \eq{Euler-hypoplast}, we use also by-part integration in time and
the Green formula also for $\DTFe$ in the evolution rule \eq{Euler2=hypoplast}
tested by a smooth $\widetilde\SS$ with $\widetilde\SS(T)=0$
together $\vv{\cdot}\nn=0$, we obtain
\begin{align*}\nonumber
&\int_0^T\!\!\!\int_\varOmega
\DTFe{:}\widetilde\SS\,\d x\d t
=\int_0^T\!\!\!\int_\varOmega\bigg(\pdt\Fe+(\vv{\cdot}\Nabla)\Fe%-(\nabla\vv)\Fe
\bigg){:}\widetilde\SS\,\d x\d t
=\int_0^T\!\!\!\int_\varGamma(\vv{\cdot}\nn)(\Fe{:}\widetilde\SS)\,\d S\d t
\\[-.2em]&\qquad\qquad\qquad-
\int_0^T\!\!\!\int_\varOmega\Fe{:}\pdt{\widetilde\SS}
+({\rm div}\,\vv)\Fe{:}\widetilde\SS
+\Fe{:}((\vv{\cdot}\Nabla)\widetilde\SS)
\,\d x\d t
-\!\int_\varOmega\!\Fezero{:}\widetilde\SS(0)\,\d x\,.
\end{align*}
Actually, $\Fe{:}\pdt{}\widetilde\SS+\Fe{:}((\vv{\cdot}\Nabla)\widetilde\SS)$
can be written ``elegantly'' as $\DT{\widetilde\SS}$ but it combines the testing
$\widetilde\SS$ with the solution $\vv$ and we will better not use such a
``too compact'' form.

\begin{definition}[Weak solutions to \eq{Euler-hypoplast}]\label{def-hypoplast}
A triple $(\vv,\Fe,\Lp)\in L^2(I;W^{2,p}(\varOmega;\R^d))$ $\times$
\linebreak $L^\infty(I;L^r(\varOmega;\R^{d\times d}))$ $\times$
$L^\infty(I;W^{1,q}(\varOmega;\R_{\rm dev}^{d\times d}))$ will be called a weak solution
to the system \eq{Euler-hypoplast} with the boundary conditions \eq{Euler-hypoplast-BC}
and  the initial condition $\Fe|_{t=0}^{}=\Fezero$ if $\vv{\cdot}\nn=0$,
$\det\Fe>0$ a.e.\ such that $\varrho=\varrho_0/\det\Fe\in L^\infty(I{\times}\varOmega)$,
and the integral identities 
\begin{subequations}\label{Euler-weak}\begin{align}
&\nonumber
\int_0^T\!\!\!\!\int_\varOmega
\big(\varphi'(\Fe)\Fe^\top\!+\xi'(\ee(\vv))\big){:}\Nabla\widetilde\vv
-\varphi(\Fe)({\rm div}\,\widetilde\vv)
+\NU|\nabla\ee(\vv)|^{p-2}\nabla\ee(\vv)\Vdots\Nabla\ee(\widetilde\vv)\,\d x\d t
\\[-.4em]&\hspace{11em}
+\int_0^T\!\!\!\!\int_\varGamma\kappa\vv{\cdot}\widetilde\vv\,\d S\d t
=\!\int_0^T\!\!\!\!\int_\varOmega\varrho\GRAVITY{\cdot}\widetilde\vv\,\d x\d t
+\int_0^T\!\!\!\!\int_\varGamma\ff{\cdot}\widetilde\vv\,\d S\d t
\\[-2.5em]\nonumber
\intertext{and}\nonumber
&\int_0^T\!\!\!\int_\varOmega\bigg(\Fe{:}\pdt{\widetilde\SS}+\Big(({\rm div}\,\vv)
\Fe{+}(\Nabla\vv)\Fe{-}\Fe\Lp\Big){:}\widetilde\SS
\\[-.8em]&\hspace{16em}
+\Fe{:}((\vv{\cdot}\Nabla)\widetilde\SS)\bigg)\,\d x\d t
=-\!\!\int_\varOmega\!\Fezero{:}\widetilde\SS(0)\,\d x
\label{Euler2=hypoplast-weak}
\intertext{hold for any $\widetilde\vv$ and $\widetilde\SS$ smooth with
$\widetilde\vv{\cdot}\nn=0$, $\widetilde\vv(T)=0$, and $\widetilde\SS(T)={\bm0}$,
and also the variational inequality}
&\int_0^T\!\!\!\int_\varOmega\zeta(\widetildeLp)
+\frac{\MU}q|\Nabla\widetildeLp|^q
-\Fee^{\top}\SS{:}(\widetildeLp{-}\Lp)\,\d x\d t
\label{Euler3=hypoplast-weak}
\ge\int_0^T\!\!\!\int_\varOmega\zeta(\Lp)+\frac{\MU}q|\Nabla\Lp|^q\,\d x\d t
\end{align}
\end{subequations}
holds for any $\widetildeLp\in L^\infty(I;W^{1,q}(\varOmega;\R_{\rm dev}^{d\times d}))$.
\end{definition}

Before stating the main analytical result, let us summarize the data
qualification (for some $\epsilon>0$):
\begin{subequations}\label{Euler-ass}\begin{align}
&\varOmega\ \text{ a smooth bounded domain of $\R^d$, }\ d=2,3
\\&\label{Euler-ass-phi}
\varphi\in C^1({\rm GL}^+(d)),\ \ \varphi(F)\ge\epsilon/(\det F)^\varkappa,
\\&\xi\in C^1(\R_{\rm sym}^{d\times d})\ \text{ convex}\,,\ \xi(0)=0\,,\
\mbox{$\sup_{\R_{\rm sym}^{d\times d}}$}
|\xi'(\cdot)|/(1+|\cdot|^{p-1})<\infty\,,
\label{Euler-ass-xi}
\\&\nonumber\zeta:\R^{d\times d}\to[0,+\infty]\
\text{satisfy \eq{ass-inelastic-dissip},
$\ \forall L\in\R_{\rm dev}^{d\times d}\!:\
\lambda\mapsto\zeta(\lambda L)$ is differentiable at $\lambda=1$,}
\\&\hspace*{1em}%\nonumber
\mbox{$\exists q_0\ge1:\ \ \inf_{\R_{\rm dev}^{d\times d}{\setminus}\{0\}}$}
\zeta(\cdot)/
%\ge\epsilon
|\cdot|^{q_0}>0\ \ \text{ and }\ \
\mbox{$\sup_{\R_{\rm dev}^{d\times d}}$}|\pl\xi(\cdot)|/(1+|\cdot|^{q-1})<\infty\,,
\label{Euler-ass-zeta}
\\&\kappa\in L^\infty(\varGamma)\,,\ \ {\rm ess\hspace*{.05em}inf}\,\kappa>0\,,
\ \ \NU>0\,,\ \MU>0\,,
\\&%\GRAVITY\in L^{p'}(I;L^1(\varOmega))+L^1(I;L^2(\varOmega))
\GRAVITY\in L^{2\varkappa/(\varkappa-2)}
%^{\varkappa p/(\varkappa p-\varkappa-p)}
(I;L^{\varkappa'}(\varOmega;\R^d))\,,
\ \ \ff\in L^2(I{\times}\varGamma;\R^d)\,,
%\ \text{ with }\ff{\cdot}\nn=0,
\label{Euler-ass-f-g}
\\&\Fezero\in W^{1,r}(\varOmega;\R^{d\times d})\ \ \ \text{ with }\ \ \
{\rm ess\hspace*{.05em}\inf}_{\varOmega}^{}\det\Fezero>0,
\label{Euler-ass-Fe0}
\\&\varrho_0\in L^\infty(\varOmega)\cap W^{1,1}(\varOmega)\ \ \text{ with }\ \ 
{\rm ess\hspace*{.05em}\inf}_{\varOmega}^{}\varrho_0>0\,.
\end{align}\end{subequations}
where ${\rm GL}^+(d)=\{F\in\R^{d\times d};\ \det F>0\}$ denotes the orientation-preserving
general linear group. As we admit $q_0=1$ in \eqref{Euler-ass-zeta}, the ``essential part'' of the dissipation potential \eqref{dissip-pot-grad} can be degree-1 homogeneous,
which would lead to a rate-independent plasticity like e.g.\
\cite{GraSte17FPPP,MaiMie09GERI,Miel03EFME,MieRou16RIEF} although all other
dissipative mechanisms stay rate dependent. 

\begin{proposition}[Existence and regularity of weak solutions]\label{prop-Euler-plast}
Let $\min(p,q)>d$ and the assumptions \eq{Euler-ass} hold
for $\varkappa>rd/(r{-}d)$  with some $r>d$. Then:\\
\Item{(i)}{there exist a weak solution according Definition~\ref{def-hypoplast}
such that also 
$\pdt{}{\Fe}\in L^2(I;L^r(\varOmega;\R^{d\times d}))$ and
$\Nabla\Fe\in L^\infty(I;L^r(\varOmega;\R^{d\times d\times d}))$.
Moreover, it conserves energy in the sense that the energy 
dissipation  balance \eq{energy-hypoplast} integrated over time interval
$[0,t]$ with the initial condition $\Fe|_{t=0}^{}=\Fezero$ holds.}
\Item{(ii)}{If also $\Fpzero\in L^s(\varOmega;\R^{d\times d})$ with some
$s>1$ and $\det\Fpzero=1$ a.e.\ on $\varOmega$, then
the corresponding plastic distortion $\Fp$ reconstructed
as a unique weak solution to \eq{plastic-strain-reconstructed}
belongs to $L^\infty(I;L^s(\varOmega;\R^{d\times d}))$ and 
$\det\Fp=1$ a.e.\ on $I{\times}\varOmega$.}
\Item{(iii)}{If also $\Fpzero\in W^{1,s}(\varOmega;\R_{\rm dev}^{d\times d})$ with
some $s>1$, then the plastic distortion $\Fp$ belongs to
$L^\infty(I;W^{1,s}(\varOmega;\R_{\rm dev}^{d\times d}))\,\cap\,W^{1,p}(I;L^s(\varOmega;\R_{\rm dev}^{d\times d}))$ and the deformation gradient
$\FF=\Fe\Fp\in L^\infty(I;W^{1,\min(s,s^*r/(s^*+r))}(\varOmega;\R^{d\times d}))$ with 
$\pdt{}\FF\in
L^2(I;L^{s^*r/(s^*+r)}(\varOmega;\R^{d\times d}))+L^p(I;L^s(\varOmega;\R^{d\times d}))$,
and ${\rm ess\hspace*{.2em}inf}_{I\times\varOmega}^{}\,\det\FF>0$.}
\end{proposition}

\begin{proof}
For clarity, we will divide the proof into nine steps.

\medskip\noindent{\it Step 1: a regularization and discretization}.
Let us first make formally the  a~priori  estimates which follow
from the energetics \eq{energy-hypoplast} when use the assumptions
\eq{Euler-ass} for $\varkappa>rd/(r{-}d)$  with some $r>d$
and $\min(p,q)>d$ and the Healey-Kr\"omer Lemma.

The only difficult term is $\varrho\GRAVITY{\cdot}\vv$ on the right-hand side of
\eq{energy-hypoplast},
which can be estimated by H\"older's and Young's inequalities as
\begin{align}\nonumber
&\!\!\int_\varOmega\varrho\GRAVITY{\cdot}\vv\,\d x
=\int_\varOmega\frac{\varrho_0}{\det\Fe}\GRAVITY{\cdot}\vv\,\d x
\le\bigg\|\frac{\varrho_0}{\det\Fe}\bigg\|_{L^{\varkappa}(\varOmega)}
\big\|\vv\big\|_{L^\infty(\varOmega;\R^d)}\big\|\GRAVITY\big\|_{L^{\varkappa'}(\varOmega;\R^d)}
\\&\le C_\epsilon
\bigg(1+\Big\|\frac{\varrho_0}{\det\Fe}\Big\|_{L^{\varkappa}(\varOmega)}^{\varkappa}\!
+\big\|\GRAVITY\big\|_{L^{\varkappa'}(\varOmega;\R^d)}^{2\varkappa/(\varkappa-2)}\bigg)
+\epsilon\big\|\Nabla\ee(\vv)\big\|_{L^p(\varOmega;\R^{d\times d\times d})}^p\!+
\epsilon\big\|\vv|_\varGamma^{}\big\|_{L^2(\varGamma;\R^d)}^2\,,
\label{est-of-rho.f.v}\end{align}
where we used $\|\vv\|_{L^\infty(\varOmega;\R^d)}\le
C(\|\Nabla\ee(\vv)\|_{L^p(\varOmega;\R^{d\times d\times d})}^{}
+\|\vv|_\varGamma^{}\|_{L^2(\varGamma;\R^d)}^{})$;
here a Korn-Poincar\'e  inequality with the Navier boundary condition
for $\vv$ is  exploited; for even a stronger variant exploiting only an
$L^2$-norm of the deviatoric part of $\ee(\vv)$ with $L^1$-norm of 
the instead of the $L^p$-norm of
$\Nabla\ee(\vv)$ cf.\ \cite[Theorems~10.16--10.17]{FeiNov00SLTV}.
The  term $\|\varrho_0/\det\Fe\|_{L^{\varkappa}(\varOmega)}^{\varkappa}$
in \eq{est-of-rho.f.v} can thus be treated by the Gronwall inequality
relying on the blow-up assumption $\varphi(F)\ge\epsilon/(\det F)^\varkappa$ 
in \eq{Euler-ass-phi}. Of course, we choose $\epsilon>0$ sufficiently small
so that the last and the penultimate terms in \eq{est-of-rho.f.v} can be
absorbed in the left-hand side of the energy balance.
 From \eq{energy-hypoplast}, we thus obtain 
\begin{subequations}\label{Euler-quasistatic-est-formal}
\begin{align}\label{Euler-plast-quasistatic-est-formal1}
&\|\vv\|_{%L^\infty(I;L^2(\varOmega;\R^d))\cap
L^2(I;W^{2,p}(\varOmega;\R^d))}^{}\le C\ \ \ \text{ with }\ \ \
\|\Nabla\ee(\vv)\|_{L^p(I{\times}\varOmega;\R^{d\times d\times d})}^{}\le C\,,
\\&\label{Euler-plast-quasistatic-est-formal3}
\|\Lp\|_{L^{q_0}(I{\times}\varOmega;\R_{\rm dev}^{d\times d})}^{}\le C\,
\ \ \,\text{ with }\ \ \
\|\Nabla\Lp\|_{L^q(I{\times}\varOmega;\R^{d\times d\times d})}^{}\le C\,.
\intertext{and $\int_\varOmega\varphi(\Fe(t))\,\d x$ bounded uniformly in
time. The exponent $2$ in the $L^2$-estimate %\eq{Euler-est-Galerkin-v}
\eq{Euler-plast-quasistatic-est-formal1} is due to the linearity of the Navier
boundary condition  and thus due to the only quadratic growth of the term
$\kappa|\vv|^2$ in \eq{energy-hypoplast}.
The former estimate in \eq{Euler-plast-quasistatic-est-formal1} together
with the qualification \eq{Euler-ass-Fe0} of the initial condition $\Fezero$
can then be exploited for the estimation as \eq{test-FF+} and
\eq{test-Delta-r} below to obtain $W^{1,r}$-regularity of $\Fe$. 
The usage of Lemma~\ref{lem-Healey-Kromer} gives $1/\det\Fe(t)$ bounded
in $L^\infty(\varOmega)$. Altogether,  for some $\varepsilon>0$, we obtain}
&\label{Euler-plast-quasistatic-est-formal2}
\|\Fe\|_{L^\infty(I;W^{1,r}(\varOmega;\R^{d\times d}))}\le C\ \ \text{ with }\ \ %\det\Fe\ge
{\rm ess\hspace*{.05em}inf}_{I\times\varOmega}^{}\det\Fe>\varepsilon\,.%\ \ \text{ and }\ \ 
\intertext{As $r>d$, \eq{Euler-plast-quasistatic-est-formal2} implies also}
&|\Fe|<\frac1\varepsilon,\ \ \ \ 
\varphi(\Fe)<\frac1\varepsilon,\ \ \text{ and }\ \ |\varphi'(\Fe)|<\frac1\varepsilon\ \ \text{ a.e.\ on }\ I{\times}\varOmega\,;
\label{Euler-plast-quasistatic-est-formal4}
\intertext{without loss of generality, we may take $\varepsilon>0$
small enough so that  both  \eq{Euler-plast-quasistatic-est-formal2} and
 \eq{Euler-plast-quasistatic-est-formal4} hold with the same $\varepsilon$.
Thus,  $|\varphi'(\Fe)\Fe^\top|\le1/\varepsilon^2$ and 
$|\Fe^\top\varphi'(\Fe)|\le1/\varepsilon^2$. 
Moreover, from \eq{Euler3=hypoplast} by comparison  realizing that the
Eshelby stress on the right-hand side is bounded on $I{\times}\varOmega$,
we can even improve
%\eq{Euler-est-Galerkin-L}
 the time integrability of  \eq{Euler-plast-quasistatic-est-formal3} as}
&\label{Euler-est-Galerkin-L-better}
\|\Lp\|_{L^\infty(I;W^{1,q}(\varOmega;\R_{\rm dev}^{d\times d}))}^{}\le C\,.
\end{align}\end{subequations}

Then, taking this $\varepsilon>0$ from (\ref{Euler-quasistatic-est-formal}c,{d}),
we make a regularization of the right-hand side of the momentum
equation and the Cauchy and the Eshelby stresses and, for %some $\delta>0$
$k\in\N$, a parabolic regularization of the evolution equation \eq{Euler2=hypoplast}
for $\Fe$. Altogether, we devise the regularized system as 
\begin{subequations}\label{Euler-plast-quasistatic-regul}\begin{align}
\nonumber
&{\rm div}{\bm\varSigma}+\varrho\GRAVITY={\bm0}\ \ \ \ \ \text{ with }\ \
\varrho=\frac{\varrho_0}{\max(\det\Fe,\varepsilon)}\ \ \text{and }\ 
     \ \ \ {\bm\varSigma}=\TT-{\rm div}\,\mathfrak{H},
    \\\nonumber
    &\hspace*{7.2em}\text{where }\
    \TT=\frac{\SS\Fee^\top}{1+(|\SS\Fee^{\top}|\,{-}1/\varepsilon^2)^+\!\!}
    +\frac{\varphi(\Fe)\bbI}{1+(\varphi(\Fe){-}1/\varepsilon)^+\!\!}+\xi'(\ee(\vv))
    \\
    &\hspace*{7.2em}
 \text{with }\ \
\SS=\varphi'(\Fe)\ \ \text{ and }\ \ 
\mathfrak{H}=\NU|\nabla\ee(\vv)|^{p-2}\nabla\ee(\vv)\,,
\label{Euler-plast-quasistatic-regul1}
\\\label{Euler-plast-quasistatic-regul2}
&\DTFe=(\nabla\vv)\Fe-\Fe\Lp+k^{-1}{\rm div}(|\Nabla\Fe|^{r-2}\Nabla\Fe)\,,
\\\label{Euler-plast-quasistatic-regul3}
&\pl\zeta(\Lp)-%\MU\Delta\Lp
{\rm div}(\MU|\Nabla\Lp|^{q-2}\Nabla\Lp)\ni %-
\frac{\Fee^{\top}\SS}{1+(|\Fee^{\top}\SS|{-}1/\varepsilon^2)^+\!\!}%\,.
\end{align}\end{subequations}
with $(\cdot)^+=\max(\cdot,0)$. The boundary conditions \eq{Euler-hypoplast-BC}
must now be complemented by some boundary condition for the regularizing term in 
 \eq{Euler-plast-quasistatic-regul2}, say $\nn{\cdot}\Nabla\Fe={\bf0}$. 
Of course, we will be interested in weak solutions to
\eq{Euler-plast-quasistatic-regul} with these boundary conditions and the
initial condition for $\Fe$. The corresponding weak formulation a'la
Definition~\ref{def-hypoplast} is quite straightforward and we will not 
explicitly write it, also because it is obvious from its Galerkin
version \eq{Euler-weak-Galerkin} below.  The philosophy of the
regularization \eq{Euler-plast-quasistatic-regul} is that the estimation of
\eq{Euler-plast-quasistatic-regul1} and (\ref{Euler-plast-quasistatic-regul}b,c)
decouples and simultaneously the a~priori estimates are the same as the
formal estimates \eq{Euler-quasistatic-est-formal} and, when taking
$\varepsilon>0$ small to comply with (\ref{Euler-quasistatic-est-formal}c,d),
the $\varepsilon$-regularization becomes eventually inactive, cf.\ Step~7 below.
Moreover, the parabolic regularization of the flow rule
\eq{Euler-plast-quasistatic-regul2} can be suppressed, cf.\ Step 6. 

Then we make a  conformal  Galerkin approximation of
\eq{Euler-plast-quasistatic-regul1} by using a nested finite-dimensional
subspaces $\{V_k\}_{k\in\N}$ whose union is dense in
$W^{2, p}(\varOmega;\R^d)$; note that they are indexed by the
same $k\in\N$ as used in \eq{Euler-plast-quasistatic-regul2}. Separately, we
make a Galerkin approximation of \eq{Euler-plast-quasistatic-regul2} and
\eq{Euler-plast-quasistatic-regul3} by using other nested finite-dimensional
subspaces $\{W_{ l}\}_{ l\in\N}$  whose union is dense in
$W^{1,\max(q,r)}(\varOmega;\R^{d\times d})$, using another index $l\in\N$.
 Also, the trace-free functions from $\{W_{l}\}_{l\in\N}$ are dense in
in the space $\{\bm{L}\in W^{1,q}(\varOmega;\R^{d\times d});\ {\rm tr}(\bm{L})=0\}$.
Without loss of generality, we may assume $\vv\in V_1$ and  $\Fezero\in W_1$. 

The approximate solution of the regularized system
\eq{Euler-plast-quasistatic-regul} will be denoted by
$(\vv_{kl},\Fekl,\Lpkl):I\to V_k\times W_l\times W_l$.  Specifically, such
a triple should satisfy the following integral identities
\begin{subequations}\label{Euler-weak-Galerkin}\begin{align}
&\nonumber
\int_0^T\!\!\!\!\int_\varOmega\bigg(\frac{\varphi'(\Fekl)\Fekl^\top}
{1+(|\varphi'(\Fekl)\Fekl^\top|{-}1/\varepsilon^2)^+\!}
\!+\xi'(\ee(\vv_{kl}))\bigg){:}\Nabla\widetilde\vv
+\frac{\varphi(\Fekl){\rm div}\,\widetilde\vv}{1+(\varphi(\Fekl){-}1/\varepsilon)^+}
\\[-.1em]&\hspace{5em}\nonumber
+\NU|\nabla\ee(\vv_{kl})|^{p-2}\nabla\ee(\vv_{kl})\Vdots\Nabla\ee(\widetilde\vv)%\bigg)
\,\d x\d t+\int_0^T\!\!\!\!\int_\varGamma\kappa\vv_{kl}{\cdot}\widetilde\vv\,\d S\d t
\\[-.4em]&\hspace{11.5em}
=\!
\int_0^T\!\!\!\!\int_\varOmega\frac{\varrho_0\GRAVITY}{\max(\det\Fekl,\varepsilon)}{\cdot}\widetilde\vv\,\d x\d t
+\int_0^T\!\!\!\!\int_\varGamma\ff{\cdot}\widetilde\vv\,\d S\d t\,,
\intertext{and}
\nonumber&\int_0^T\!\!\!\int_\varOmega\bigg(\Fekl{:}%\DT{\widetilde\SS\!\!\!}
\pdt{\widetilde\SS}+\Big(({\rm div}\,\vv_{kl})
\Fekl{+}(\Nabla\vv_{kl})\Fekl{-}\Fekl\Lpkl\Big){:}\widetilde\SS
\\[-.3em]&\hspace{2em}
+\Fekl{:}((\vv_{kl}{\cdot}\Nabla)\widetilde\SS)
-\frac1k|\Nabla\Fekl|^{r-2}\Nabla\Fekl\Vdots\Nabla\widetilde\SS
\bigg)\,\d x\d t
=-\!\!\int_\varOmega\!\Fezero{:}\widetilde\SS(0)\,\d x
\label{Euler2=hypoplast-weak-Galerkin}
\intertext{for any $\widetilde\vv\in L^\infty(I;V_k)$ and
 $\widetilde\SS\in L^\infty(I;W_l)$ with $\widetilde\vv{\cdot}\nn=0$,
$\widetilde\vv(T)=0$, and $\widetilde\SS(T)={\bm0}$,
and also the variational inequality}
\nonumber
&\int_0^T\!\!\!\int_\varOmega\zeta(\widetildeLp)
+\frac{\MU}q|\Nabla\widetildeLp|^q
-%\Fekl^{\top}\SS
\frac{\Fekl^\top\varphi'(\Fekl)}{1+(|\Fekl^\top\varphi'(\Fekl)|{-}1/\varepsilon^2)^+}
{:}(\widetildeLp{-}\Lpkl)\,\d x\d t
%\\[-.4em]&\hspace{10em}
\\[-.2em]&\hspace{17em}
\ge
\int_0^T\!\!\!\int_\varOmega\zeta(\Lpkl)+\frac{\MU}q|\Nabla\Lpkl|^q\,\d x\d t
\label{Euler3=hypoplast-weak-Galerkin}
\end{align}
\end{subequations}
should hold for any $\widetildeLp\in L^\infty(I;W_l)$, ${\rm tr}(\widetildeLp)=0$ a.e.\
on $I{\times}\Omega$.

Existence of this solution is based on the theory of systems of
ordinary differential equations first locally in time, and then
by successive prolongation on the whole time interval based on 
the $L^\infty$-estimates below.

\medskip\noindent{\it Step 2: first  a~priori  estimates}.
The basic test of the Galerkin approximation of
\eq{Euler-plast-quasistatic-regul} can be done 
by $(\vv_{kl},\Fekl,\Lpkl)$.
In particular, for \eq{Euler-plast-quasistatic-regul2}, we use
the estimate
\begin{align}\nonumber
&\frac{\d}{\d t}\int_\varOmega\frac12|\Fekl|^2\,\d x
+\frac1k\int_\varOmega|\Nabla\Fekl|^r\,\d x
\\[-.4em]&\nonumber\quad\ \ 
\le\int_\varOmega\Big((\nabla\vv_{kl})\Fekl-(\vv_{kl}{\cdot}\nabla)\Fekl-\Fekl\Lpkl\Big){:}\Fekl\,\d x
\\[-.4em]&\nonumber\qquad\ \ 
=\int_\varOmega\!(\nabla\vv_{kl})\Fekl{:}\Fekl
+\frac{{\rm div}\,\vv_{kl}}2|\Fekl|^2
-\Lpkl{:}(\Fekl^\top\Fekl)\,\d x
\\[-.4em]&\qquad\qquad\ \ \le
\bigg(\frac32\|\nabla\vv_{kl}\|_{L^\infty(\varOmega;\R^{d\times d})}^{}
+\|\Lpkl\|_{L^\infty(\varOmega;\R^{d\times d})}^{}\bigg)
\|\Fekl\|_{L^2(\varOmega;\R^{d\times d})}^2\,;
\label{test-FF}\end{align}
%for any $\delta>0$;
here we used also the calculus (for $\FF=\Fekl$ and $\vv=\vv_{kl}$)
\begin{align}\nonumber
  \int_\varOmega(\vv{\cdot}\nabla)\FF{:}\FF\,\d x
  &=\int_\varGamma|\FF|^2(\vv{\cdot}\nn)\,\d S
  \\[-.4em]&\quad
  -\int_\varOmega\!\FF{:}(\vv{\cdot}\nabla)\FF+({\rm div}\,\vv)|\FF|^2\,\d x
=-\frac12\int_\varOmega({\rm div}\,\vv)|\FF|^2\,\d x
\label{calcul-FF}\end{align}
together with the boundary condition $\vv{\cdot}\nn=0$. By the %discrete
Gronwall inequality  exploiting the first left-hand-side term
which does not contain the factor $1/k$, we obtain the estimate
\begin{align}
\big\|\Fekl\big\|_{L^\infty(I;L^2(\varOmega;\R^{d\times d}))}^{}\le %C
\|\Fezero\|_{L^2(\varOmega;\R^{d\times d}))}^{}
{\rm e}^{\|\nabla\vv_{kl}\|_{L^1(I;L^\infty(\varOmega;\R^{d\times d}))}
+\|\Lpkl\|_{L^1(I;L^\infty(\varOmega;\R^{d\times d}))}}\,.
\label{est-plast-Euler}\end{align}
Thus, by this test, we obtain
\begin{subequations}\label{Euler-quasistatic-est1}
\begin{align}\label{Euler-quasistatic-est1-1}
&\|\vv_{kl}\|_{L^{2}(I;W^{2,p}(\varOmega;\R^d))}\le C\,,
\\&\label{Euler-quasistatic-est1-2}
\|\Fekl\|_{L^\infty(I;L^2(\varOmega;\R^{d\times d}))}\le C\ \ \text{ with }\ \ \|\Nabla\Fekl\|_{L^r(I{\times}\varOmega;\R^{d\times d\times d})}\le C\sqrt[r]{k}\,,
\\&\label{Euler-quasistatic-est1-3}
\|\Lpkl\|_{L^\infty(I;W^{1,^q}(\varOmega;\R_{\rm dev}^{d\times d}))}\le C\,.
\end{align}\end{subequations}
Particularly, let us note that the regularized force $\varrho_0\GRAVITY/\max(\det\Fekl,\varepsilon)$ is  a~priori  bounded in $L^{2\varkappa/(\varkappa-2)}(I;L^{\varkappa'}(\varOmega;\R^d))$, cf.\ \eq{Euler-ass-f-g},
and that the constant $C$ in \eq{Euler-quasistatic-est1-1}
depends on $\varepsilon$ but not on $k,\ l,$ %and $\delta$
because also the  conservative ({elastic}) part of the regularized
Cauchy stress in \eq{Euler-plast-quasistatic-regul1}
and also the Eshelby stress in \eq{Euler-plast-quasistatic-regul3}
are bounded independently of $\Fe$, and therefore each of the
equations in \eq{Euler-plast-quasistatic-regul} can be estimated separately.

\medskip\noindent{\it Step 3: second  a~priori  estimates}.
In Step 1, we could also estimate $\pdt{}\Fekl+(\vv_{kl}{\cdot}\Fekl
-(\nabla\vv_{kl})\Fekl+\Fekl\Lpkl$ by comparison from
\eq{Euler-quasistatic-est1-2} in $L^{r'}(I;W^{1,r}(\varOmega;\R^{d\times d})^*)$, but
it would not be enough for \eq{Euler-quasistatic-est3-1} below and thus for
making the test in Step~5 legitimate. To get a better estimate, we can also test the
Galerkin approximation of \eq{Euler-plast-quasistatic-regul2} by $\pdt{}\Fekl$.
By H\"older and Young inequalities, we can estimate
\begin{align}\nonumber
&\int_\varOmega\bigg|\pdt{\Fekl}\bigg|^2\,\d x
+\frac%\delta
1{rk}\frac{\d}{\d t}\int_\varOmega|\Nabla\Fekl|^r\,\d x
\\[-.1em]&\nonumber\ 
\le\int_\varOmega\Big((\nabla\vv_{kl})\Fekl-(\vv_{kl}{\cdot}\nabla)\Fekl-\Fekl\Lpkl\Big){:}\pdt{\Fekl}\,\d x
\\[-.4em]&\nonumber\qquad 
\le\Big(\|\nabla\vv_{kl}\|_{L^\infty(\varOmega;\R^{d\times d})}^{}
+\|\Lpkl\|_{L^\infty(\varOmega;\R^{d\times d})}^{}\Big)^2
\|\Fekl\|_{L^2(\varOmega;\R^{d\times d})}^2
\\[-.2em]&\qquad\qquad\ 
+C_r\|\vv_{kl}\|_{L^\infty(\varOmega;\R^d)}^2
\Big(%|\varOmega|
1+\|\Nabla\Fekl\|_{L^r(\varOmega;\R^{d\times d})}^r\Big)
+\frac12\bigg\|\pdt{\Fekl}\bigg\|_{L^2(\varOmega;\R^{d\times d})}^2
\label{test-FF+}
\end{align}
with some $C_r\in\R$; here we used that surely $r>2$. Using the already obtained
estimates \eq{Euler-quasistatic-est1} and the Gronwall inequality, we obtain
\begin{align}\label{Euler-quasistatic-est2}
&
\Big\|\pdt{\Fekl}\Big\|_{L^2(I{\times}\varOmega;\R^{d\times d}))}\le C{\rm e}^{k}/k\ \ \text{ and }\ \
\|\Nabla\Fekl\|_{L^\infty(I;L^r(\varOmega;\R^{d\times d\times d}))}\le C
{\rm e}^{k/r}\,.
\end{align}%\end{subequations}
 Note that here the Gronwall inequality uses not the first but the second
left-hand-side term which contains the factor $1/k$ so that both
estimates in \eq{Euler-quasistatic-est2} are $k$-dependent. 

\medskip\noindent{\it Step 4: limit passage with $l\to\infty$}.
Now, by the Banach selection principle, we extract a subsequence and some
$(\vv_k,\FEk,\Lpk):I\to V_k\times W^{1,r}(\varOmega;\R^{d\times d})\times W^{1,q}(\varOmega;\R_{\rm dev}^{d\times d})$ such that
\begin{subequations}\label{Euler-plast-weak}
\begin{align}
\label{Euler-weak-v}
&\!\!\vv_{kl}\to\vv_k&&\text{weakly* in $\ L^2(I;W^{2,p}(\varOmega;\R^d))$,}\!\!&&
\\
&\!\!\Fekl\to\FEk\!\!\!&&\text{weakly* in $\ 
L^\infty(I;W^{1,r}(\varOmega;\R^{d\times d}))\,\cap\,H^1(I;L^2(\varOmega;\R^{d\times d}))$,}\!\!
\\
&\!\!\Lpkl\to\Lpk\!\!\!&&\text{weakly* in $\ L^\infty(I;W^{1,q}(\varOmega;\R_{\rm dev}^{d\times d}))$}
\,.
\intertext{By the Aubin-Lions theorem, we have also
$\Fekl\to\FEk$ strongly in $L^c(I{\times}\varOmega;\R^{d\times d})$ with any
$c<\infty$; recall that $r>d$. Thus, by the continuity of the corresponding Nemytski\u{\i}
(or here simply superposition) mappings, also the regularized Eshelby stress
converges}
&\nonumber
\frac{\Fekl^\top\varphi'(\Fekl)}{1+(|\Fekl^\top\varphi'(\Fekl)|{-}1/\varepsilon^2)^+}
\to
\frac{\FEk^\top\varphi'(\FEk)}{1+(|\FEk^\top\varphi'(\FEk)|{-}1/\varepsilon^2)^+}
\hspace*{-20em}&&\hspace*{20em}\text{strongly in $L^c(I{\times}\varOmega;\R^{d\times d})$.}
\intertext{As the right-hand side of the discretized plasticity-rate inclusion
\eq{Euler-plast-quasistatic-regul3} converge strongly, by the uniform monotonicity of its
left-hand side, we have also} 
&\!\!\Lpkl\to\Lpk\!\!\!&&\text{strongly in
$\ L^c(I;W^{1,q}(\varOmega;\R_{\rm dev}^{d\times d}))$ with any $c<\infty$}
\,.
\end{align}\end{subequations}
By the mentioned continuity of the corresponding Nemytski\u{\i} mappings, we have also
\begin{align*}
\varrho_{kl}=\frac{\varrho_0}{\max(\det\Fekl,\varepsilon)}
\to\frac{\varrho_0}{\max(\det\FEk,\varepsilon)}=\varrho_k\ \
\ \text{ strongly in $L^c(I{\times}\varOmega)$ with any $c<\infty$}
\end{align*}
and similarly the regularized elastic part of the Cauchy stress in
\eq{Euler-plast-quasistatic-regul1} converges strongly:
\begin{align}\nonumber
&\TT_{\varepsilon,kl}=
\frac{\varphi'(\Fekl)\Fekl}{1+(|\varphi'(\Fekl)\Fekl|{-}1/\varepsilon^2)^+\!}
+\frac{\varphi(\Fekl)\bbI}{1+(\varphi(\Fekl){-}1/\varepsilon)^+}
\\&\ \ \to
\frac{\varphi'(\FEk)\FEk}{1{+}(|\varphi'(\FEk)\FEk|{-}1/\varepsilon^2)^+\!}
+\frac{\varphi(\FEk)\bbI}{1{+}(\varphi(\FEk){-}1/\varepsilon)^+}
=\TT_{\varepsilon,k}
\hspace*{1em}\text{strongly in $L^c(I{\times}\varOmega;\R^{d\times d})$.}
\label{Euler-Cauchy-stong-conv}\end{align}

For the limit passage in the Galerkin approximation of 
\eq{Euler-plast-quasistatic-regul1} for $l\to\infty$ (which still
will remain discretized as $k$ is considered fixed in this step), we need
the Minty trick or, having here the strong monotonicity of the
hyperstress term, we use just strong convergence of $\Nabla\ee(\vv_{kl})$.
In fact, as we do not consider any acceleration and thus do not have
$\pdt{}\vv_{kl}$ under control in the quasistatic case, we will anyhow need
strong convergence of $\Nabla\vv_{kl}$ later
in Step~6 for the convective term $(\vv{\cdot}\Nabla)\Fe$. So, 
using the Galerkin approximation of the momentum equation tested by $\vv_{kl}-\vv_k$,
we can estimate
\begin{align}\nonumber
  &\NU c_p%\lim_{h\to0}
  \|\nabla\ee(\vv_{kl}{-}
\vv_k)\|_{L^p(I{\times}\varOmega;\R^{d\times d\times d})}^p\le%\lim_{h\to0}
\int_0^T\!\!\!\int_\varOmega\!\Big(\big(\xi'(\ee(\vv_{kl})){-}
\xi'(\ee(\vv_k))\big){:}\ee(\vv_{kl}{-}\vv_k)
  \\[-.4em]&\hspace*{4em}\nonumber
 +\NU\big(|\nabla\ee(\vv_{kl})|^{p-2}\nabla\ee(\vv_{kl})
-|\nabla\ee(\vv_k)|^{p-2}\nabla\ee(\vv_k)\big)\Vdots
  \nabla\ee(\vv_{kl}{-}\vv_k)\Big)\,\d x\d t
  \\[-.5em]&=\nonumber
 \bigg(\int_0^T\!\!\!\int_\varOmega\Big(\varrho_{kl}\GRAVITY{\cdot}(\vv_{kl}{-}\vv_k)
 -\NU\big(|\nabla\ee(\vv_k)|^{p-2}\nabla\ee(\vv_k)\big)\Vdots
 \nabla\ee(\vv_{kl}{-}\vv_k)
 \\[-.5em]&\hspace*{3em}
-\big(\TT_{\varepsilon,kl}+\xi'(\ee(\vv_k))\big){:}\Nabla(\vv_{kl}{-}\vv_k)\Big)\,\d x\d t
+\int_0^T\!\!\!\int_\varGamma\ff{\cdot}(\vv_{kl}{-}\vv_k)\,\d S\d t
\bigg)\to0\,
\label{strong-hyper+}\end{align}
with some $c_{p}>0$ related to the inequality
$c_{p}|G-\widetilde G|^p\le(|G|^{p-2}G-|\widetilde G|^{p-2}\widetilde G)\Vdots(G-\widetilde G)$ holding for $p\ge2$.
We also use \eq{Euler-Cauchy-stong-conv} and that $\nabla(\vv_{kl}{-}\vv_k)\to0$ weakly
in $L^p (I{\times}\varOmega;\R^{d\times d})$, so that
$\int_0^T\!\int_\varOmega\TT_{\varepsilon,kl}{:}\nabla(\vv_{kl}{-}\vv_k)\,\d x\d t\to0$,
and also the growth assumption \eq{Euler-ass-xi} which ensures
$\xi'(\ee(\vv_k))\in L^p(I{\times}\varOmega;\R^{d\times d})$.
Thus we obtain the desired strong convergence of
$\nabla\ee(\vv_{kl})$ in $L^p(I{\times}\varOmega;\R^{d\times d\times d})$.

The limit passage in the quasilinear parabolic evolution equation
\eq{Euler-plast-quasistatic-regul2} in its Galerkin approximation
 \eq{Euler2=hypoplast-weak-Galerkin}  is very standard when realizing
that $(\Nabla\vv_{kl})\Fekl-\Fekl\Lpkl$ converges to $(\Nabla\vv_k)\FEk-\FEk\Lpk$
strongly in $L^p(I;L^{q^*}(\varOmega;\R^{d\times d}))$ 
while $(\vv_{kl}{\cdot}\Nabla)\Fekl\to(\vv_k{\cdot}\Nabla)\FEk$
only weakly* in $L^\infty(L^r(\varOmega;\R^{d\times d}))$ but this is enough when tested
by $\widetildeFekl-\Fekl\to0$ strongly in $L^1(I;L^{r'}(\varOmega;\R^{d\times d}))$
with some approximation $\widetildeFekl\in L^1(I;W_l)$ strongly converging to
$\FEk$ for $l\to\infty$; the strong convergence $\Fekl\to\FEk$ is due to the
Aubin-Lions theorem and the estimates \eq{Euler-quasistatic-est2}; here
we rely on that surely $1/r+1/r^*<1$.

The limit passage in the Galerkin approximation of \eq{Euler-plast-quasistatic-regul3}
written as the variational inequality \eq{Euler3=hypoplast-weak-Galerkin}  can be made
easily by a weak convergence and by the weak lower semicontinuity of the functional
$\Lp\mapsto\int_0^T\!\int_\varOmega\zeta(\Lp)+\MU|\Nabla\Lp|^q/q\,\d x\d t$.

\medskip\noindent{\it Step 5: third  a~priori  estimates}. Since now
\begin{align}
\Big\|\pdt{\FEk}+
(\vv_k{\cdot}\nabla)\FEk-
(\nabla\vv_k)\FEk+\FEk\Lpk\Big\|_{L^2(I{\times}\varOmega;\R^{d\times d})}\!\le C{\rm e}^{k/r}
\,,\end{align}
by comparison we also obtain
\begin{align}\label{Euler-quasistatic-est3-1}
&\big\|{\rm div}(|\Nabla\FEk|^{r-2}\Nabla\FEk^{})\big\|_{L^2(I{\times}\varOmega;\R^{d\times d})}^{}\le
 k C{\rm e}^{k/r}\,.
\end{align}
Although this estimate blows up when $k\to\infty$, 
we have now at least the information 
that ${\rm div}(|\Nabla\FEk|^{r-2}\Nabla\FEk)\in
L^2(I{\times}\varOmega;\R^{d\times d})$.  It is now important
that  we have
%\eq{Euler-quasistatic-regul2}
\eq{Euler-plast-quasistatic-regul2} continuous, i.e.\ non-discretized.
Therefore, we can legitimately
use ${\rm div}(|\Nabla\FEk|^{r-2}\Nabla\FEk)$ as a test.
Since $\min(p,q)>d$, we have $\min(p,q)^{-1}+(r^*)^{-1}+(r')^{-1}\le1$,
and thus by the H\"older and Young inequalities, we can estimate
\begin{align}\nonumber
&\frac{\d}{\d t}
\int_\varOmega\frac1r|\nabla\FEk|^r\,\d x
+%\delta
\frac1k\int_\varOmega|{\rm div}(|\Nabla\FEk|^{r-2}\Nabla\FEk)|^2\,\d x
\\[-.0em]\nonumber&
=\int_\varOmega\nabla\Big((\vv_k{\cdot}\nabla)\FEk
  -(\nabla\vv_k)\FEk-\FEk\Lpk\Big)\Vdots\big(|\nabla\FEk|^{r-2}\nabla\FEk\big)\,\d x
\\[-.1em]\nonumber&=
\int_\varOmega|\nabla\FEk|^{r-2}(\nabla\FEk{\otimes}\nabla\FEk){:}\ee(\vv_k)-\frac1r|\nabla\FEk|^r{\rm div}\,\vv_k
\\[-.4em]&\hspace*{2em}\nonumber
-\Big((\nabla\vv_k)\nabla\FEk+(\nabla^2\vv_k)\FEk+\Nabla\FEk\Lpk
+\FEk\Nabla\Lpk\Big)\Vdots\big(|\nabla\FEk|^{r-2}\nabla\FEk\big)\,\d x
\\&\nonumber
\le C_r\Big(\|\nabla\vv_k\|_{L^\infty(\varOmega;\R^{d\times d})}^{}
+\|\Lpk\|_{L^\infty(\varOmega;\R^{d\times d})}^{}\Big)
\|\nabla\FEk\|_{L^r(\varOmega;\R^{d\times d\times d})}^r\!
\\[-.2em]&\nonumber\ \ \ +%\delta
C_r\Big(\|\nabla^2\vv_k\|_{L^p(\varOmega;\R^{d\times d\times d})}^{}\!
+\|\Nabla\Lpk\|_{L^q(\varOmega;\R^{d\times d\times d})}^{}\Big)
\|\FEk\|_{L^{r^*}(\varOmega;\R^{d\times d})}^{}%\Big(1{+}
\|\nabla\FEk\|_{L^{r}(\varOmega;\R^{d\times d\times d})}^{r-1}
\\&\le\nonumber
C_r\Big(\|\nabla\vv_k\|_{L^\infty(\varOmega;\R^{d\times d})}^{}
+\|\Lpk\|_{L^\infty(\varOmega;\R^{d\times d})}^{}\Big)
\|\nabla\FEk\|_{L^r(\varOmega;\R^{d\times d\times d})}^r\!
\\[-.2em]&\nonumber\ \ \ +%\delta
C_rN\Big(\|\nabla^2\vv_k\|_{L^p(\varOmega;\R^{d\times d\times d})}^{}\!\!
+\|\Lpk\|_{L^2(\varOmega;\R^{d\times d})}^{}\Big)
\|\FEk\|_{L^2(\varOmega;\R^{d\times d})}
\big(1{+}\|\nabla\FEk\|_{L^{r}(\varOmega;\R^{d\times d\times d})}^r\big)
\\&\ \ \ +
C_rN\Big(\|\nabla^2\vv_k\|_{L^p(\varOmega;\R^{d\times d\times d})}^{}\!
+\|\Nabla\Lpk\|_{L^q(\varOmega;\R^{d\times d\times d})}^{}\Big)
\|\nabla\FEk\|_{L^{r}(\varOmega;\R^{d\times d\times d})}^r\,,
\label{test-Delta-r}
\end{align}
where we used $\min(p,q)>d$ also for the embedding of $\Nabla\vv_k$
and $\Lpk$ into $L^\infty(\varOmega;\R^{d\times d})$ and where we further used
the calculus (to be used for $\FF=\FEk$)
\begin{align}\nonumber
&\int_\varOmega\nabla\big((\vv{\cdot}\nabla)\FF
  \big){:}|\nabla\FF|^{r-2}\nabla\FF\,\d x
 \\[-.4em]&\hspace{2em}\nonumber=\int_\varOmega|\nabla\FF|^{r-2}(\nabla\FF{\otimes}\nabla\FF){:}\ee(\vv)  
+(\vv{\cdot}\nabla)\nabla\FF\Vdots|\nabla\FF|^{r-2}\nabla\FF\,\d x
\\[-.1em]&\hspace{2em}\nonumber=\int_\varGamma|\nabla\FF|^r\vv{\cdot}\nn\,d S
+\int_\varOmega\bigg(|\nabla\FF|^{r-2}(\nabla\FF{\otimes}\nabla\FF){:}\ee(\vv)
\\[-.8em]&\hspace{12em}\nonumber
-({\rm div}\,\vv)|\nabla\FF|^r-(r{-}1)|\nabla\FF|^{r-2}\nabla\FF\Vdots
(\vv{\cdot}\nabla)\nabla\FF\bigg)\,\d x
\\[-.4em]&\hspace{2em}\nonumber
=\int_\varGamma\frac{|\nabla\FF|^r\!\!}r\ \vv{\cdot}\nn\,d S
+\int_\varOmega|\nabla\FF|^{r-2}(\nabla\FF{\otimes}\nabla\FF){:}\ee(\vv)-({\rm div}\,\vv)\frac{|\nabla\FF|^r\!\!}r\ \d x\,.
%\label{test-Delta+}
\end{align}
Here $\nabla\FF{\otimes}\nabla\FF$ denoted the
symmetric matrix $[\nabla\FF{\otimes}\nabla\FF]_{ij}^{}=
\sum_{k,l=1}^d\frac{\partial}{\partial\xx_i}{\FF}_{kl}^{}
\frac{\partial}{\partial\xx_j}{\FF}_{kl}^{}$.
Again, the boundary integral vanishes in \eq{test-Delta-r} if
$\vv{\cdot}\nn=0$. For the last inequality in \eq{test-Delta-r}, we %will use
 have used  $\|\FEk\|_{L^{r^*}(\varOmega;\R^{d\times d})}^{}\le
N(\|\FEk\|_{L^2(\varOmega;\R^{d\times d})}^{}+\|\Nabla\FEk\|_{L^{r}(\varOmega;\R^{ d\times d\times d})}^{})$
where $N$ is the norm of the embedding $W^{1,r}(\varOmega)\subset
L^{r*}(\varOmega)$ if $W^{1,r}(\varOmega)$ is endowed with the norm
$\|\cdot\|_{L^2(\varOmega)}+\|\nabla\cdot\|_{L^r(\varOmega;\R^d)}$.

Thus one can apply  the  Gronwall inequality to \eq{test-Delta-r} and the estimates 
\eq{Euler-quasistatic-est1-2} and \eq{Euler-quasistatic-est2} can be strengthened.
Specifically, using the already obtained estimates \eqref{est-plast-Euler} and having
assumed $\Fezero\in W^{1,r}(\varOmega;\R^{d\times d})$, one obtains the estimates
\begin{subequations}\label{Euler-quasistatic-est3}
\begin{align}\label{Euler-quasistatic-est3-1+}
&\|\Nabla\FEk\|_{L^\infty(I;L^r(\varOmega;\R^{d\times d\times d}))}\le C\ \ \text{ and }
\\&\label{Euler-quasistatic-est3-2}
\|{\rm div}(|\Nabla\FEk|^{r-2}\Nabla\FEk^{})\|_{L^2(I{\times}\varOmega;\R^{d\times d})}\le C\sqrt k\,.
\intertext{Besides,  although  the former estimate in
\eq{Euler-quasistatic-est2} on $\pdt{}\FEk$
is  not  inherited also on the limit, we have by comparison
from $\pdt{}\FEk=(\nabla\vv_k)\FEk-(\vv_k{\cdot}\nabla)\FEk
+k^{-1}{\rm div}(|\Nabla\FEk|^{r-2}\Nabla\FEk^{})$ in its Galerkin approximation
\eq{Euler2=hypoplast-weak-Galerkin} a weaker estimate}
&\Big\|\pdt{\FEk}\Big\|_{L^2(I; L^2(\varOmega;\R^{d\times d})+W^{1,r}_K(\varOmega;\R^{d\times d})^*)}\le C\ \
\text{ for }\ k\ge K\,,
\label{Euler-quasistatic-est3-2+}\end{align}\end{subequations}
 where $W^{1,r}_K(\Omega;\R^{d\times d})^*$ is considered endowed with the seminorm
$$
|\cdot|_K^{}=\sup_{\|F\|_{W^{1,r}(\Omega;\R^{d\times d})}^{}\le1,\ F\in W_K}\int_\varOmega\nabla\cdot\Vdots\nabla F\,\d x\,.
$$
It is important that $C$ in \eq{Euler-quasistatic-est3-2+}
can be taken independent of $K\in\N$. 

\medskip\noindent{\it Step 6: limit passage with %$\delta\to0$
$k\to\infty$}.
We use the Banach selection principle as in Step~4 now also taking 
\eq{Euler-quasistatic-est3-1+} into account instead of 
the latter estimate in \eq{Euler-quasistatic-est2} which was not uniform
in $k$. Thus, for a subsequence and some $(\vv,\Fe,\Lp)$, we have 
\begin{subequations}\label{Euler-plast-weak+}
\begin{align}
&\!\!\vv_k\to\vv&&\text{weakly* in $\ L^2(I;W^{2,p}(\varOmega;\R^d))$,}\!\!&&
\\
&\!\!\FEk\to\Fe\!\!\!&&\text{weakly* in $\ L^\infty(I;W^{1,r}(\varOmega;\R^{d\times d})))$,
}\!\!
\label{Euler-plast-weak+2}\\
&\!\!\Lpk\to\Lp\!\!\!&&\text{strongly in
$\ L^c(I;W^{1,q}(\varOmega;\R_{\rm dev}^{d\times d}))$ with any $c<\infty$}
\,.
\end{align}\end{subequations}
 Moreover, exploiting \eq{Euler-plast-weak+2} together with the estimate
\eq{Euler-quasistatic-est3-2+}, by the Aubin-Lions theorem generalized
for time derivatives controlled in Hausdorff locally convex spaces
\cite[Lemma~7.7]{Roub13NPDE} we obtain also
$\FEk\to\Fe$ strongly in $L^c(I{\times}\varOmega;\R^{d\times d})$ for 
any $1\le c<+\infty$ to be used analogously as we did in Step~4. 

The momentum equation \eq{Euler-plast-quasistatic-regul1} (still regularized
by $\varepsilon$ and discretised) and the plastic-rate
inclusion \eq{Euler-plast-quasistatic-regul3} are to be treated like
in Step~4; in fact, \eq{strong-hyper+} is to be slightly modified by using
some approximation $\widetilde\vv_k$ of the limit $\vv$
valued in the Galerkin finite-dimensional space so that
$\vv_k{-}\widetilde\vv_k$ is a legitimate test for the Galerkin approximation
of the momentum equation \eq{Euler-plast-quasistatic-regul1}
and such that $\Nabla\ee(\widetilde\vv_k)\to\Nabla\ee(\vv)$
strongly in $L^p(I{\times}\varOmega;\R^{d\times d\times d})$. Due to
\eq{Euler-quasistatic-est3-2}, we have
$%\delta
k^{-1}{\rm div}(|\Nabla\FF_k|^{r-2}\Nabla\FF_k^{})=\mathscr{O}(\!\sqrt{1/k})\to0$
in $L^2(I{\times}\varOmega;\R^{d\times d})$ and thus this regularizing
term in the elastic-strain evolution equation
\eq{Euler-plast-quasistatic-regul2} disappears in the limit. The rest is a
linear equation in terms of $\Fe$, while its coefficients $\vv_k$,
$\Nabla\vv_k$, and $\Lpk$ converge strongly.
Altogether, we showed that  $(\vv,\Fe,\Lp)$ is a weak solution of a problem like
\eq{Euler-plast-quasistatic-regul} but regularized only by
$\varepsilon>0$, i.e.\ the last term in \eq{Euler-plast-quasistatic-regul2}
is omitted.

\medskip\noindent{\it Step 7: the original problem}.
Let us note that  the limit $\Fe$ lives in 
$L^\infty(I;W^{1,r}(\varOmega))\,\cap\,H^1(I;W^{1,r}(\varOmega)^*)$  and this space
is embedded into $C(I{\times}\barOmega)$ if $r>d$. Therefore $\Fe$ and its determinant
evolve continuously in time, being valued respectively in $C(\barOmega;\R^{d\times d})$
and $C(\barOmega)$. Let us recall that the initial condition $\Fezero$ complies with
the bounds (\ref{Euler-quasistatic-est-formal}c,{d}) and we used this $\Fezero$
also for the $\varepsilon$-regularized system.
Therefore $\Fe$ satisfies these bounds not only at $t=0$ but also at least
for small times. Yet, it means that the $\varepsilon$-regularization is nonactive
and  $(\vv,\Fe,\Lp)$ solves, at least for a small time, the original
nonregularized system for which the  a~priori  bounds \eq{Euler-quasistatic-est-formal}
hold. Here we used Lemma~\ref{lem-Healey-Kromer}. By the continuation argument, we may
see that the $\varepsilon$-regularization remains therefore inactive within the whole
evolution of $(\vv,\Fe,\Lp)$ on the whole time interval $I$.

\medskip\noindent{\it Step 8: energy balance}.
It is now important that the tests and then all the subsequent calculations
leading to \eq{energy-hypoplast} integrated over a current
time interval $[0,t]$ are really legitimate.

Since $\Nabla\ee(\vv)\in L^{p}(I{\times}\varOmega;\R^{d\times d\times d})$, we have 
${\rm div}^2(\NU|\Nabla\ee(\vv)|^{p-2}\Nabla\ee(\vv))\in
L^{p'}(I;W^{2,p}(\varOmega;\R^d)^*)$ in
duality with $\vv$. Also ${\rm div}\xi'(\ee(\vv))\in
L^{p'}(I;W^{1,p^*}(\varOmega;\R^d)^*)$ is in duality with $\vv$ due to the
growth condition \eq{Euler-ass-xi} and
$\varrho_0\GRAVITY/\det\Fe+{\rm div}(\varphi'(\Fe)\Fe^\top+\varphi(\Fe)\bbI)$ is even
better. Further, by comparison, $\pdt{}\Fe=(\Nabla\vv)\Fe
-(\vv{\cdot}\Nabla)\Fe-\Fe\Lp\in L^p(I;L^r(\varOmega;\R^{d\times d}))+
L^q(I{\times}\varOmega;\R^{d\times d})\subset L^{\min(p,q)}(I;L^{\min(r,q)}(\varOmega;\R^{d\times d}))$.
Therefore, it is surely in duality with the Piola stress
$\SS=\varphi'(\Fe)\in L^\infty(I{\times}\varOmega;\R^{d\times d})$.
Also, by comparison from \eq{Euler3=hypoplast},
${\rm div}(\MU|\Nabla\Lp|^{q-2}\Nabla\Lp)\in \pl\zeta(\Lp)+\Fee^{\top}\SS$ is a
bounded set in $L^{q'}(I;L^\infty(\varOmega;\R_{\rm dev}^{d\times d}))$, cf.\ the
growth condition in \eq{Euler-ass-zeta}. Therefore, it  is  in duality
with $\Lp\in L^q(I;W^{1,q}(\varOmega;\R_{\rm dev}^{d\times d}))$. Here we note that,
due to \eq{Euler-ass-zeta}, $\pl\zeta$ is possibly multivalued but, due to
\eq{Euler-ass-zeta},  the plastic dissipation rate  %only at $0$ so that
$\pl\zeta(\Lp){:}\Lp$ is, in fact, always a single-valued function
in $L^1(I{\times}\varOmega)$.

Therefore, the calculations
\eq{Euler-hypoplast-test-momentum}--\eq{Euler-hypoplast-test-momentum++}
are legitimate.

\medskip\noindent{\it Step 9: additional information -- plastic distortion}.
The corresponding plastic distortion $\Fp$ satisfies the evolution rule
$\DTFp=\Lp\Fp$, cf.\ \eq{plastic-strain-reconstructed}.
Then, for the $W^{1,s}$-estimate, it suffices to apply
the same procedure as we did for \eq{Euler2=hypoplast} modified and even simplified
since there is no term like $(\Nabla\vv)\Fp$, i.e.\ we regularize as the linear
transport-evolution equation $\DTFp=\Lp\Fp+k{\rm div}(|\Nabla\Fp|^{s-2}|\Nabla\Fp)$
%instead of $\DTFe=(\Nabla\vv)\Fe-\Fe\Lp$.
and use the calculus like \eq{test-FF}, \eq{test-FF+}, and \eq{test-Delta-r}.
Actually, for the mere $L^s$-estimate with $s=2$, it
suffices to use only the first estimate \eq{test-FF} while, for $L^s$-estimate
with $s\ne2$, the estimate \eq{test-FF} is to be modified for a test by
$|\Fp|^{s-2}\Fp$ in the spirit of \eq{test-Delta-r}. Moreover,
$\Lp\Fp\in L^\infty(I;L^{s^*\!}(\varOmega;\R^{d\times d}))$
and $(\vv{\cdot}\Nabla)\Fp\in L^p(I;L^s(\varOmega;\R^{d\times d}))$ so that 
we have $\pdt{}\Fp=\Lp\Fp-(\vv{\cdot}\Nabla)\Fp\in
L^p(I;L^s(\varOmega;\R^{d\times d}))$.

The same arguments can be applied to the { evolution-and-}transport  equation  
of $\det\Fp$, cf.\ 
\eq{D/DT-of-determinant}. Due to \eq{ass-inelastic-dissip}, ${\rm tr}\Lp=0$
so that \eq{D/DT-of-determinant} reduces to 
\begin{align}\label{D/DT-of-determinant+}
\DT{\overline{\det\Fp}}=0\,.
\end{align}
If $\det\Fpzero$ is constant, then \eq{D/DT-of-determinant+} reduces
to $\pdt{}(\det\Fp)=0$, so that $\det\Fp$ stays equal to this constant
during the whole evolution. In particular it holds for $\det\Fpzero=1$.

Eventually, for the deformation gradient
$\FF=\Fe\Fp\in L^\infty(I;L^{s^*\!}(\varOmega;\R^{d\times d}))$, we have also
$\Nabla\FF=\Nabla\Fe\Fp+\Fe\Nabla\Fp\in
L^\infty(I;L^{\min(s,s^*r/(s^*+r))}(\varOmega;\R^{d\times d\times d}))$ and
$\pdt{}\FF=(\pdt{}\Fe)\Fp+\Fe\pdt{}\Fp\in
L^2(I;L^{s^*r/(s^*+r)}(\varOmega;\R^{d\times d}))+L^p(I;L^s(\varOmega;\R^{d\times d}))$.
As $\det\Fp=1$ and $\det\Fe$ stays away from 0, the same holds for
$\det\FF=\det\Fe$.
\end{proof}

\begin{remark}[{\sl Classical solutions}]\upshape
In fact, we proved that  all the  terms in the transport-evolution
equations \eq{Euler2=hypoplast} and, under the assumptions of
Proposition~\ref{prop-Euler-plast}(iii), also \eq{plastic-strain-reconstructed}
are surely in $L^1(I{\times}\varOmega;\R^{d\times d})$. Thus these equations are
satisfied even a.e.\ on $I{\times}\varOmega$.  If $\varphi$ is
twice continuously differentiable,
${\rm div}(\varphi'(\Fe)\Fe^{\top}\!+\varphi(\Fe)\bbI)\in
L^\infty(I;L^r(\varOmega;\R^d))$ due to the regularity of
$\Nabla\Fe\in L^\infty(I;L^r(\varOmega;\R^{d\times d}))$. If also $\zeta$ is
twice continuously differentiable, then also ${\rm div}(\zeta'(\ee(\vv))\in
L^p(L^\infty(\varOmega;\R^d))$. Then, by comparison,
${\rm div}^2(\NU|\Nabla\ee(\vv)|^{p-2}\Nabla\ee(\vv))\in
L^p(I;L^r(\varOmega;\R^d))$ and therefore also the momentum
equation \eq{Euler1=hypoplast} holds a.e.\ on $I{\times}\varOmega$.
If the plastic-distortion-rate inclusion \eq{Euler3=hypoplast} is understood
on the linear subspace ${\rm tr}\,\Lp={\bm0}$, then also
${\rm div}(\MU|\Nabla\Lp|^{q-2}\Nabla\Lp)\in
L^{q'}(I{\times}\varOmega;\R^{d\times d})$ and also the inclusion
\eq{Euler3=hypoplast} holds a.e.\ on $I{\times}\varOmega$.
 This   is more than the weak formulation \eq{Euler-weak}.
 Recovery of the boundary conditions a.e.\ on $I{\times}\varGamma$ would
need still more regularity, however.

\end{remark}

\begin{remark}[{\sl Uniqueness}]\upshape
For a given $\vv$ and $\Lp$, the weak solution of the transport-evolution
equations \eq{Euler2=hypoplast} is unique. The highest-order terms of the
momentum equation \eq{Euler1=hypoplast} and the plastic-distortion-rate
inclusion \eq{Euler3=hypoplast} are strictly monotone but, anyhow, the
uniqueness of
a weak solution to the whole system \eq{Euler-hypoplast} seems problematic.
The troublesome attribute is that the conservative part of the Cauchy stress
$\varphi'(\Fe)\Fe^{\top}\!+\varphi(\Fe)\bbI$ as well as the Eshelby
stress $\Fe^{\top}\!\varphi'(\Fe)$ are highly nonmonotone.
\end{remark}

\bigskip

{\small

\baselineskip=12pt

\noindent{\it Acknowledgments.}
The author is very thankful for extremely valuable and
inspiring discussions and  comments to the manuscript to 
Giuseppe Tomassetti. Also valuable discussions about transport
equations with Sebastian Schwarzacher  and about hypoplasticity models
with Yannis F.\ Dafalias  are warmly acknowledged.
 Careful reading and many valuable suggestions by two anonymous
referees are thankfully acknowledged, too. 
Also the supports from the M\v SMT \v CR (Ministry of Education of the
Czech Republic) project CZ.02.1.01/0.0/0.0/15-003/0000493
and the institutional support RVO:61388998 (\v CR) are acknowledged.

%\bibliographystyle{plain}
%\bibliography{tr-plast-Eulerian.bib}

}

\end{document}